\documentclass[ejs]{imsart}

\usepackage{amsthm,amsmath,natbib}
\usepackage{amssymb,mathrsfs}
\usepackage{graphicx}

\RequirePackage[dvips]{hyperref}

\RequirePackage{hypernat}

\pubyear{2008}

\startlocaldefs
\newtheorem{theo}{Theorem}

\newtheorem{prop}{Proposition}
\newtheorem{lemm}{Lemma}
\newtheorem{defi}{Definition}
\newcommand{\pen}{\text{\upshape{pen}}}
\newcommand{\cdi}{\ensuremath{c_0}}
\newcommand{\argmin}[1]{\underset{#1}{\text{argmin}}}
\newdimen\AAdi%
\newbox\AAbo%
\def\AArm{\fam0 }
\def\AAk#1#2{\setbox\AAbo=\hbox{#2}\AAdi=\wd\AAbo\kern#1\AAdi{}}%
\def\BBone{{\AArm 1\AAk{-.8}{I}I}}%

\renewcommand{\thetable}{\thesection.\arabic{table}}

\endlocaldefs

\begin{document}

\begin{frontmatter}

\title{Detecting change-points in a discrete distribution via model selection}
\runtitle{Detecting change-points in a discrete distribution}


\begin{aug}
\author{\fnms{Nathalie} \snm{Akakpo}\ead[label=e1]{nathalie.akakpo@math.u-psud.fr}}
\address{Laboratoire de Probabilit\'es et Statistiques\\
B\^at. 425, Universit\'e Paris-Sud, 91405 Orsay Cedex, France\\\printead{e1}}
\end{aug}

\runauthor{N. Akakpo}

\begin{abstract}
This paper is concerned with the detection of multiple change-points in the joint distribution of independent categorical variables. The procedures introduced rely on model selection and are based on a penalized least-squares criterion. Their performance is assessed from a nonasymptotic point of view. Using a special collection of models, a preliminary estimator is built. According to an existing model selection theorem, it satisfies an oracle-type inequality. Moreover, thanks to an approximation result demonstrated in this paper, it is also proved to be adaptive in the minimax sense. In order to eliminate some irrelevant change-points selected by that first estimator, a two-stage procedure is proposed, that also enjoys some adaptivity property. Besides, the first estimator can be computed with a complexity only linear in the size of the data. A heuristic method allows to implement the second procedure quite satisfactorily with the same computational complexity.
\end{abstract}





\end{frontmatter}

\section{Introduction}

Let $Y_1, Y_2,\ldots, Y_n$ be independent random variables taking value in the finite set $\{1,\ldots, r\}$, where $r$ is an integer and $r\geq 2$, and let $s$ be the joint distribution of $(Y_1, Y_2,\ldots, Y_n)$. Assume that $\{1,\ldots,n\}$ can be partitioned into intervals such that all the $Y_i$'s with indices $i$ in a same interval follow the same law. Then $s$ is said to have change-points located at the beginning of each interval, 1 excluded. In this paper, our aim is to detect change-points in $s$, using no a priori information on their number. A typical example of application is given by the DNA segmentation problem, for which the review~\cite{BraunMuller} by Braun and M\"uller may serve as an introduction. The $n$-uple $(Y_1, Y_2,\ldots, Y_n)$ provides indeed a model for the successive bases along a DNA sequence of length $n$, when coding the set of bases \{Adenine, Cytosine, Guanine, Thymine\} by $\{1,\ldots, 4\}$ for instance. Thus, beyond the theoretical properties of the statistical procedures, a special attention must be paid to their computational complexity, due to the length of sequences such as DNA ones.

Several methods based on a penalized criterion, with a penalty typically increasing with the number of change-points, have been proposed for the statistical problem under consideration.
Braun, Braun and M\"uller present in~\cite{BraunBraunMuller} such a procedure, based on a penalized quasi-deviance criterion, and prove consistency results for the estimation of the change-points and the true number of change-points. Nevertheless, the computational complexity of their estimator, though reduced by using dynamic programming, is quite costly, with $\mathcal{O}(n^3)$ computations, or $\mathcal{O}(n^2D_{max})$ if an upper-bound $D_{max}$ is imposed on the number of change-points.
Lebarbier and N\'ed\'elec also study penalized criteria in~\cite{LebarbierNédélec}, one based on least-squares, the other on maximum likelihood. Their procedures are based on the model selection principle developed by Birg\'e and Massart in various papers, such as~\cite{BirgéMassart}. Thus they adopt a wholly different point of view from that of Braun et al.: the estimators studied in~\cite{LebarbierNédélec} are nonparametric and are proved to satisfy a nonasymptotic oracle-type inequality, for an adequate choice of the penalty. But, when considering all possible configurations of change-points, these procedures suffer from the same computational complexity as that of Braun et al.
In view of significantly reducing the computational time, the CART-based procedure proposed by Gey and Lebarbier in a Gaussian regression framework (cf.~\cite{GeyLebarbier}) can be adapted to the framework considered here, as illustrated in~\cite{Lebarbier}, Chapter 7. In the best case, the number of computations falls down to only $\mathcal{O}(n\ln(n))$. Unfortunately, apart from the the oracle-type inequality given in~\cite{LebarbierNédélec}, theoretical properties of that hybrid procedure seem difficult to establish. Adopting the same approach as in~\cite{LebarbierNédélec} or~\cite{BirgéMassart}, Durot, Lebarbier and Tocquet propose in~\cite{DLT} quite a general framework for estimating $s$ relying on a penalized least-squares criterion, where the choice of the penalty is supported by an oracle-type inequality. As a particular case, Durot et al. recover one of the change-point detection methods proposed in~\cite{LebarbierNédélec}. They complete the study of its performance with an improved oracle-type inequality and an adaptivity result in the minimax sense.
Let us also mention some other methods, not based on penalized criteria, that enjoy some interesting computational complexity. They are not supported however by theoretical results.
Fu and Curnow propose in~\cite{FuCurnow} an estimator based on maximum likelihood, imposing a constraint on the minimal lengths of the segments to prevent overfitting. According to~\cite{Csuros}, it can be implemented with a computational complexity only linear in the size of the data. Szpankowski, Szpankowski and Ren study in~\cite{SzpanSzpanRen} a procedure inspired from Information Theory. It also has a linear complexity, that  results from the splitting of the sequence into blocks of a prescribed length.

Following the work presented in~\cite{Lebarbier},~\cite{LebarbierNédélec} and~\cite{DLT}, we propose in this paper two statistical procedures based on a penalized least-squares criterion, using the same model selection principle. Each estimator we build is piecewise constant on a partition of $\{1,\ldots,n\}$. If the distribution $s$ is piecewise constant, then the partition associated with the estimator allows to estimate its change-points. We first study an estimator based on a special collection of models in correspondence with the partitions of $\{1,\ldots,n\}$ into dyadic intervals only. That collection of models satisfies two important properties.
On the one hand, it has been chosen for its potential qualities of approximation. They have been suggested by a theorem due to DeVore and Yu (cf.~\cite{DeVoreYu}) about the approximation of functions in Besov spaces by piecewise polynomials. Adapting their proof to our framework, we prove that our collection of models has indeed good approximation qualities with respect to Besov bodies, some discrete analogues of balls in a Besov space defined in this article.
On the other hand, the number of models per dimension is much lower for that collection than for the analogous one associated with all the partitions of $\{1,\ldots,n\}$ into intervals, also called exhaustive collection in~\cite{LebarbierNédélec} and~\cite{DLT}. So no extra logarithmic factor appears in the oracle-type inequality satisfied by our first estimator.
The conjunction of both properties allows to prove an adaptivity result in the minimax sense over Besov bodies. Notice that, because of those two interesting properties, a similar collection of models has lately been used by Birg\'e (cf.~\cite{BirgéTest} and~\cite{BirgéPoisson}) and Baraud and Birg\'e (cf.~\cite{BaraudBirgé}) for estimation by model selection in various statistical frameworks.
About our first procedure, we must underline that considering such a reduced collection of partitions also happens to reduce the computational complexity to the so wanted linear complexity.
It should also be noted that the hypothesis that $s$ is piecewise constant is not used to derive any result. Therefore, whatever $s$, that first procedure still provides an interesting estimator of $s$.
For the detection of change-points, if it does detect some relevant ones, it also selects some less significant ones, due to the nature of the selected partition. That's why we propose the following hybrid procedure. A preliminary stage consists in using part of the data to select a partition into dyadic intervals with the previous procedure, that will henceforth be called preliminary procedure. During the second stage, the rest of the data is used to select, among the rougher partitions built on the previous one, the one minimizing a penalized least-squares criterion. The resulting hybrid estimator also enjoys some adaptivity property, similar to that of the first procedure, up to a $\ln(n)$ factor. Moreover, in practice, it can also be implemented quite efficiently with a linear complexity.

The paper is organized as follows. In the brief section~\ref{sec:notations}, we describe the statistical framework and introduce notation used throughout the paper. The next two sections are devoted to the theoretical study of the preliminary estimator and of the subsequent hybrid estimator. The performance of these procedures are illustrated in section~\ref{sec:simulations} through a simulation study. In particular, we discuss there the practical choice of the penalties constants. The paper ends with the proof of the approximation result needed to derive the adaptivity properties of both estimators.

\section{Framework and notation}\label{sec:notations}

        \subsection{Framework}

We observe $n$ independent random variables $Y_1,\ldots,Y_n$ defined on the same probability space $(\Omega, \mathcal{A},\mathbb{P})$ and with values in $\{1,\ldots,r\}$, where $r$ is an integer and $r\geq 2$. Moreover, we assume that $n$ is a power of 2 and write $n = 2^N$.
The distribution of the $n$-uple $(Y_1,\ldots,Y_n)$ is defined as the $r\times n$ matrix $s$ whose $i$-th column is
$$s_i = \big(\mathbb{P}(Y_i = 1)\ldots\mathbb{P}(Y_i = r)\big)^T, \text{ for } 1\leq i\leq n.$$
Observing $(Y_1,\ldots,Y_n)$ is equivalent to observing the random $r \times n$ matrix $\textit{X}$ whose $i$-th column is
$$X_i = \big(\BBone_{Y_i = 1}\ldots\BBone_{Y_i = r}\big)^T,\text{ for } 1\leq i\leq n.$$
It should be noted that the distribution $s$ to estimate is in fact the mean of $X$.

        \subsection{Notation}

Let $\mathscr{M}(r,n)$ be the set of all real matrices with $r$ rows and $n$ columns.
Given an element $t \in \mathscr{M}(r,n)$, we denote by $t^{(l)}$ its $l$-th row and by $t_i$ its $i$-th column. The space $\mathscr{M}(r,n)$ is endowed with the inner product defined by
$$\langle t,t' \rangle =\sum_{i=1}^n\sum_{l=1}^r t_i^{(l)}{t'}_i^{(l)}.$$
That product is linked with the standard inner products on $\mathbb{R}^r$ and $\mathbb{R}^n$, denoted respectively by $\langle.,.\rangle_r$ and $\langle.,.\rangle_n$, by the relations
$$\langle t,t' \rangle =\sum_{i=1}^n\langle t_i,{t'}_i\rangle_r
=\sum_{l=1}^r\langle t^{(l)},{t'}^{(l)}\rangle_n. $$
The norms induced by these products on $\mathscr{M}(r,n)$, $\mathbb{R}^r$ and $\mathbb{R}^n$ are respectively denoted  by $\|.\|$, $\|.\|_r$ and $\|.\|_n$. Another norm on $\mathscr{M}(r,n)$ appearing in this paper is
$$\|t\|_\infty:=\max \big\{|t_i^{(l)}|; 1\leq i \leq n, 1 \leq l \leq r\big\}.$$

Let us now define some subsets of $\mathscr{M}(r,n)$ of special interest. The set composed of the $r\times n$ matrices whose columns are probability distributions on $\{1, \ldots, r\}$ is denoted by $\mathscr{P}$. Given a subspace $S$ of $\mathbb{R}^n$, the notation $\mathbb{R}^r \otimes S$ stands for the linear subspace of $\mathscr{M}(r,n)$ composed of the matrices whose rows all belong to $S$.

Any vector $u$ in $\mathbb{R}^n$ is identified with the function defined from $\{1,\ldots,n\}$ into $\mathbb{R}$ and whose value in $i$ is $u_i$, for $i=1,\ldots,n$.
In particular, for any subset $I$ of $\{1,\ldots,n\}$, we will call indicator function of $I$, and denote by $\BBone_I$, the $\mathbb{R}^n$-vector whose $i$-th coordinate is equal to 1 if $i\in I$, and null otherwise.

When the distribution of $(Y_1,\ldots,Y_n)$ is given by $s$, we denote respectively by $\mathbb{P}_s$ and $\mathbb{E}_s$ the underlying probability distribution on $(\Omega^{\otimes n}, \mathcal{A}^{\otimes n})$ and the associated expectation.

Last, in the many inequalities we shall encounter, the capital letters $C, C_1, \ldots$ stand for positive constants, whose value may change from one line to another. Sometimes, their dependence on one or several parameters will be indicated. For instance, the notation $C(\alpha,p)$ means that $C$ only depends on $\alpha$ and $p$.

\section{Preliminary estimator}\label{sec:distrategy}

We study in this section a first estimator of the distribution $s$. For detecting change-points in $s$, it will be used in the next section during a preliminary stage. We begin here with the definition of that preliminary estimator: we explain the underlying model selection principle and easily justify the choice of the involved penalty thanks to~\cite{DLT}. Then, we present the main result of this paper, about the adaptivity of this estimator. It derives from an approximation result that will be proved later in the article. Last, we describe the algorithm used to compute the estimator and give its computational complexity.

        \subsection{Definition of the preliminary estimator}\label{sec:def}

Let $\mathcal{M}$ be the collection of all the partitions of $\{1,\ldots,n\}$ into dyadic intervals. In order to describe it in a more constructive way, let us introduce the complete binary tree $\mathcal{T}$ with $N+1$ levels such that:
\begin{itemize}
\item[\textbullet] the root of $\mathcal{T}$ is $(0,0)$;
\item[\textbullet] for all $j \in \{1,\ldots,N\}$, the nodes at level $j$ are indexed by the elements of the set $\Lambda(j)=\{(j,k), k=0,\ldots,2^j-1\}$;
\item[\textbullet] for all $j \in \{0,\ldots,N-1\}$ and all $k \in \{0,\ldots,2^j-1\}$, the left branch that stems from node $(j,k)$ leads to node $(j+1,2k)$, and the right one, to node $(j+1,2k+1)$.
\end{itemize}
The node set of $\mathcal{T}$ is $\mathcal{N}=\cup_{j=0}^{N}\Lambda(j)$, where $\Lambda(0)=\{(0,0)\}$. The dyadic intervals of $\{1,\ldots,n\}$ are nothing but the sets
$$I_{(j,k)}=\{k2^{N-j}+1,\ldots,(k+1)2^{N-j}\}$$
indexed by the elements of $\mathcal{N}$.
Hence we deduce a one-to-one correspondence between the partitions of $\{1,\ldots,n\}$ that belong to $\mathcal{M}$ and the subsets of $\mathcal{N}$ composed of the leaves of any complete binary tree resulting from an elagation of $\mathcal{T}$.
We consider the collection of linear spaces of the form $\mathbb{R}^r\otimes S_m$, where $m\in\mathcal{M}$ and $S_m$ is the linear subspace of $\mathbb{R}^n$ generated by the indicator functions $\{\BBone_I, I\in m\}$. In the sequel, the term "model" refers indifferently to such a subspace of $\mathscr{M}(r,n)$ or to the associated partition in $\mathcal{M}$. For all $m \in \mathcal{M}$, the least-squares estimator of $s$ in $\mathbb{R}^r \otimes S_m$ is defined by
$$\hat{s}_m = \argmin{t \in \mathbb{R}^r \otimes S_m}\|\textit{X}-t\|^2.$$

Ideally, we would like to choose a model among the collection $\mathcal{M}$ such that the risk of the associated estimator is minimal. However, determining such a model requires the knowledge of $s$. Therefore the challenge is to define a procedure $\hat{m}$, based solely on the data, that selects a model for which the risk of $\hat{s}_{\hat{m}}$ almost reaches the minimal one. In other words, the estimator $\hat{s}_{\hat{m}}$ should satisfy a so-called oracle inequality
\begin{equation}
\mathbb{E}_s \big[\|s - \hat{s}_{\hat{m}}\|^2\big] \leq C \inf_{m \in \mathcal{M}}\mathbb{E}_s \big[\|s - \hat{s}_m\|^2\big]
\label{eq:oracle}.
\notag
\end{equation}
Besides, as is usually the case, the risk of each estimator $\hat{s}_m$  breaks down into an approximation error and an estimation error roughly proportional to the dimension of the model. Indeed, for all $m \in \mathcal{M}$, the estimator $\hat{s}_m$ satisfies
\begin{equation}
\|s - \bar{s}_m\|^2 + \big(1 - \|s\|_\infty\big)D_m \leq \mathbb{E}_{s} \big[\|s - \hat{s}_m\|^2\big] \leq \|s - \bar{s}_m\|^2 + \Big(1-\frac{1}{r}\Big)D_m,
\label{eq:riskonemodel}
\notag
\end{equation}
where $\bar{s}_m$ is the orthogonal projection of $s$ on $\mathbb{R}^r \otimes S_m$ and $D_m$ is the dimension of $S_m$ (cf.~\cite{DLT}, proof of Corollary 1). Reaching the minimal risk among the estimators of the collection thus amounts to realizing the best trade-off between the approximation error and the dimension of the model, that vary in opposite ways. Therefore, we consider the data-driven procedure
$$\hat{m}=\argmin{m \in \mathcal{M}} \big\{\|\textit{X}-{\hat{s}}_m\|^2+\pen(m)\big\},$$
where $\pen : \mathcal{M}\rightarrow \mathbb{R}^+$ is called penalty function. The preliminary estimator $\tilde{s}$ of $s$ is then defined as
$$\tilde{s}={\hat{s}}_{\hat{m}}.$$

Regarding the choice of an adequate penalty, we rely on results proved in~\cite{DLT}. They provide us with the following oracle inequality, up to a quantity depending on $\|s\|_\infty$, which justifies the choice of a penalty simply linear in the dimension of the models.
\begin{prop}\label{diOracle}
Let $\pen: \mathcal{M} \rightarrow \mathbb{R}^+$ be a penalty of the form
\begin{equation}
\pen(m)=\cdi D_m,
\label{eq:linpen}
\notag
\end{equation}
where, for $m\in\mathcal{M}$, $D_m$ is the dimension of $S_m$. If \cdi\:is positive and large enough and if $\|s\|_\infty < 1$, then
\begin{equation}
\mathbb{E}_{s} \big[\|s - \tilde{s}\|^2\big] \leq C(\cdi)(1-\|s\|_\infty)^{-1}\inf_{m \in \mathcal{M}}\mathbb{E}_{s} \big[\|s - \hat{s}_m\|^2\big].
\label{eq:oracletypedi}
\end{equation}
\end{prop}
\begin{proof}
Let us introduce the subcollections of models of same dimension
\begin{equation}
\mathcal{M}_D=\{m \in \mathcal{M} \text{ s.t. } D_m = D\},\text{ for }1\leq D \leq n.
\label{eq:subcollection}
\notag
\end{equation}
We look for a penalty satisfying the hypotheses of Corollary 1 in~\cite{DLT}, otherwise said of the form
$$\pen(m)=(k_1+k_2L(D_m))D_m,$$
where $k_1$ and $k_2$ are positive constants, and $\{L(D)\}_{1\leq D\leq n}$ is a family of positive numbers, called weights, such that
\begin{equation}
\sum_{D=1}^n|\mathcal{M}_D|\exp(-DL(D))\leq 1.
\notag
\end{equation}
In fact, it is enough to require that
\begin{equation}
L(D)\geq(\ln |\mathcal{M}_D|)/D + \ln 2, \text{ for all } 1\leq D\leq n.
\label{eq:choiceweight}
\notag
\end{equation}
Since the cardinal of $\mathcal{M}_D$ is equal to the number of complete binary trees with $D$ leaves resulting from an elagation of $\mathcal{T}$, it is given by the Catalan number $D^{-1}\binom{2(D-1)}{D-1}$, and thus upper-bounded by $4^D$.
Consequently, we can set all the weights equal to a same constant. Inequality~\eqref{eq:oracletypedi} then follows from the proof of Corollary 1 in~\cite{DLT}.
\end{proof}
\noindent
From now on, we will always assume that the preliminary estimator derives from a penalty of the form $\pen(m)=\cdi D_m$, where the constant \cdi\:is positive and large enough so as to yield an oracle-type inequality.
By way of comparison, let us mention that the similar procedure based on the exhaustive collection of partitions of $\{1,\ldots,n\}$ only satisfies an oracle-type inequality such as~\eqref{eq:oracletypedi} within a $\ln(n)$ factor, owing to the greater number of models per dimension for that collection (cf.~\cite{DLT}, Proposition 1).

Last, notice that $\tilde{s}$ does not necessarily belong to $\mathscr{P}$. Nevertheless, since the vector $(1 \ldots 1)$ belongs to any $S_m$, for $m \in \mathcal{M}$, the elements in a same row of $\tilde{s}$ sum up to 1. In order to get an estimator of $s$ with values in $\mathscr{P}$, we can consider the orthogonal projection of $\tilde{s}$ on the closed convex $\mathscr{P}$, whose risk is even smaller than that of $\tilde{s}$.

        \subsection{Adaptivity of the preliminary estimator}\label{sec:adaptivity}

Though the oracle-type inequality~\eqref{eq:oracletypedi} ensures that, under a minor constraint on $s$, the estimator $\tilde{s}$ is almost as good as the best estimator in the collection $\{\hat{s}_m\}_{m \in \mathcal{M}}$, it does not provide any comparison of $\tilde{s}$ to other estimators of $s$. Therefore, we now pursue the study of $\tilde{s}$ adopting a minimax point of view. We consider a large family of subsets of $\mathscr{P}$, to be defined in the next paragraph. Let us denote by $\mathcal{S}$ some subset in that family. Our aim is to compare the maximal risk of $\tilde{s}$ when $s$ belongs to $\mathcal{S}$ to the minimax risk over $\mathcal{S}$. From  Theorem 1 in~\cite{DLT}, it easily follows that an upper-bound for the risk of $\tilde{s}$ is
\begin{equation}
\mathbb{E}_{s} \big[\|s - \tilde{s}\|^2\big] \leq C(\cdi)\inf_{1 \leq D \leq n}\Big\{\inf_{m \in \mathcal{M}_D}\|s - \bar{s}_m\|^2+ D\Big\},
\label{eq:riskupperbound}
\end{equation}
where we recall that $\mathcal{M}_D=\{m \in \mathcal{M} \text{ s.t. }\dim(S_m)=D\}$ and $\bar{s}_m$ is the orthogonal projection of $s$ on $\mathbb{R}^r \otimes S_m$.
Consequently, the approximation qualities of our family of models with respect to each subset $\mathcal{S}$ remain to be evaluated. More precisely, for each subset $\mathcal{S}$, and each dimension $D$, we shall provide upper-bounds for the approximation error $\inf_{m \in \mathcal{M}_D}\|s - \bar{s}_m\|^2$ when $s\in\mathcal{S}$.

As in~\cite{DLT}, we consider subsets of $\mathscr{P}$ whose definition is inspired from the characterization in terms of wavelet coefficients of balls in Besov spaces. In order to define them, we equip $\mathbb{R}^n$ with an orthonormal wavelet basis: the Haar basis.
\begin{defi}
Let $\Lambda = \cup_{j=-1}^{N-1}\Lambda(j)$, where $\Lambda(-1)=\{(-1,0)\}$ and $$\Lambda(j)=\{(j,k), k=0,\ldots,2^j-1\}$$
for $0\leq j\leq N-1$.
Let $\varphi : \mathbb{R} \rightarrow \{-1,1\}$ be the function with support $(0,1]$ that takes value $1$ on $(0,1/2]$ and $-1$ on $(1/2,1].$\\
If $\lambda=(-1,0)$, $\phi_\lambda$ is the vector in $\mathbb{R}^n$ whose coordinates are all equal to $1/\sqrt{n}$.\\
If $\lambda=(j,k)$, where $j \neq -1$ and $k\in\Lambda(j)$, $\phi_{\lambda}$ is the vector in $\mathbb{R}^n$ whose $i-th$ coordinate is
$$\phi_{\lambda i} = \frac{2^{j/2}}{\sqrt{n}}\varphi\bigg(2^j\frac{i}{n}-k\bigg), \text{ for } i =1,\ldots,n.$$
The functions $\{\phi_{\lambda}\}_{\lambda \in \Lambda}$ are called the Haar functions. They form an orthonormal basis of $\mathbb{R}^n$ called the Haar basis.
\end{defi}
\noindent
This basis is closely linked with the collection of partitions $\mathcal{M}$: the Haar functions from a same resolution level $j$, $ 0\leq j \leq N-1$, are indexed by the nodes at level $j$ in the tree $\mathcal{T}$ (cf. Section~\ref{sec:def}), which give the supports of these wavelets.
Besides, any element $t\in\mathscr{M}(r,n)$ can be decomposed into
\begin{equation}
t=\sum_{j=-1}^{N-1} \sum_{\lambda\in\Lambda(j)} \beta_\lambda \phi_\lambda
\notag
\end{equation}
where, for all $\lambda \in \Lambda$, $\beta_\lambda$ is the column-vector in $\mathbb{R}^r$ whose $l$-th coefficient is $\beta^{(l)}_\lambda = \langle t^{(l)},\phi_\lambda \rangle_n$, for $l = 1,\ldots, r$. So, we improperly refer to the $\beta_\lambda$'s as the wavelet coefficients of $t$. We then define Besov bodies as follows.

\begin{defi}
Let $\alpha >0$, $p >0$ and $R>0$.
The set composed of all the elements $t\in\mathscr{M}(r,n)$ such that
\begin{equation}
\Bigg(\sum_{j=0}^{N-1} 2^{jp(\alpha + 1/2 - 1/p)}\sum_{\lambda \in \Lambda(j)}\|\beta_\lambda\|_r^p\Bigg)^{1/p}\leq \sqrt{n}R
\label{eq:BesovBall},
\notag
\end{equation}
where, for $l = 1,\ldots, r$, $\beta^{(l)}_\lambda = \langle t^{(l)},\phi_\lambda \rangle_n$, is denoted by $\mathscr{B}(\alpha, p, R)$ and called a Besov body.
The set of all the elements of $\mathscr{P}$ that belong to $\mathscr{B}(\alpha, p, R)$ is denoted by $\mathscr{P}(\alpha, p, R)$.
\end{defi}
\noindent
In particular, for an element of Besov body, the size of the wavelet coefficients from a same resolution level $j$ is all the smaller as $j$ is high.
For a wide range of values of the parameter $(\alpha,p,R)$, we are able to bound the approximation errors appearing in~\eqref{eq:riskupperbound} uniformly over $\mathscr{P}(\alpha, p, R)$.
\begin{theo}\label{treeapprox}
Let $p\in (0,2]$, $\alpha > 1/p-1/2$ and $R>0$.
For all $D \in \{1,\ldots,n\}$,
\begin{equation}
\sup_{s \in \mathscr{P}(\alpha, p, R)}\inf_{m \in \mathcal{M}_D}{\|s - \bar{s}_m\|^2}
\leq C(\alpha,p)nR^2D^{-2\alpha}.
\notag
\end{equation}
\end{theo}
\noindent
That result will be proved in section~\ref{sec:proof}.

Let us now come back to our initial problem, that is comparing the performance of $\tilde{s}$ to that of any other estimator of $s$. For $\alpha>0$, $p>0$ and $R>0$, the minimax risk over $\mathscr{P}(\alpha, p, R)$ is given by
$$\mathcal{R}(\alpha, p, R)=\inf_{\hat{s}}\sup_{s \in \mathscr{P}(\alpha, p, R)}\mathbb{E}_{s} \big[\|s - \hat{s}\|^2\big]$$
where the infimum is taken over all the estimators $\hat{s}$ of $s$.
Thanks to the above approximation result, we obtain, as stated below, that, for a whole range of values of $(\alpha,p,R)$, the estimator $\tilde {s}$ reaches the minimax risk over $\mathscr{P}(\alpha, p, R)$ within a multiplicative constant. Otherwise said, $\tilde{s}$ is adaptive in the minimax sense over that range of subsets of $\mathscr{P}$.
\begin{theo}\label{adaptivity}
For all $p\in (0,2],\alpha >1/p-1/2$ and $n^{-1/2}\leq R < n^\alpha$,
\begin{equation}
\sup_{s \in \mathscr{P}(\alpha, p, R)}\mathbb{E}_{s} \big[\|s - \tilde{s}\|^2\big]
\leq C(\cdi,\alpha,p)\mathcal{R}(\alpha,p,R).
\label{eq:AdaptDyad}
\end{equation}
\end{theo}
\begin{proof}
Let us fix $p\in (0,2],\alpha >1/p-1/2$ and $n^{-1/2}\leq R < n^\alpha$. Combining Inequality~\eqref{eq:riskupperbound} and Theorem~\ref{treeapprox} leads to
\begin{equation}
\sup_{s \in \mathscr{P}(\alpha, p, R)}\mathbb{E}_{s} \big[\|s - \tilde{s}\|^2\big] \leq C(\cdi,\alpha,p)\inf_{1\leq D \leq n}\big\{nR^2D^{-2\alpha}+ D\big\}.
\label{eq:riskupperbound2}
\notag
\end{equation}
In order to realize approximately the best trade-off between the terms $nR^2D^{-2\alpha}$ and $D$, that vary in opposite ways when $D$ increases, we choose $D$ as large as possible under the constraint $D\leq nR^2D^{-2\alpha}$. Let us denote by $D^{\star}$ the largest integer $D$ such that $D\leq(nR^2)^{1/(1+2\alpha)}$. One can easily check that, given the hypotheses linking $n$ and $R$, $D^{\star}$ does belong to $\{1,\ldots,n\}$ and provides the upper-bound
\begin{equation}
\sup_{s \in \mathscr{P}(\alpha, p, R)}\mathbb{E}_{s} \big[\|s - \tilde{s}\|^2\big] \leq C(\cdi,\alpha,p)(nR^2)^{1/(2\alpha+1)}.
\label{eq:rate}
\notag
\end{equation}
The matching lower bound for the minimax risk over $\mathscr{P}(\alpha, p, R)$ has been proved in~\cite{DLT} (Theorem 3).
\end{proof}

        \subsection{Computing the preliminary estimator}\label{sec:computing}

Since the penalty only depends on the dimension of the models, we will also denote by $\pen(D)$ the penalty assigned to all models in $\mathcal{M}_D$, for $1\leq D\leq n$. A way to compute $\tilde{s}$ could rely on the equality
$$\min_{m\in\mathcal{M}} \big\{\|X-\hat{s}_m\|^2+\pen(m)\big\}
=\min_{1\leq D\leq n} \bigg\{\min_{m\in\mathcal{M}_D} \|X-\hat{s}_m\|^2+\pen(D)\bigg\}.$$
That would lead us to compute a best estimator for each dimension, before choosing one among them by taking into account the penalty term, as in~\cite{LebarbierNédélec} for the exhaustive collection of partitions or in~\cite{BraunBraunMuller}. But, even when using Bellman's algorithm, that requires polynomial time. Here, we shall see that we can avoid such a computationaly intensive way by taking advantage of the form of the penalty.

Let us express more explicitly the criterion to be minimized by $\hat{m}$. For $m\in\mathcal{M}$, we denote by $\{i_k,\ldots,i_{k+1}-1\}$, $1\leq k\leq D_m$, the dyadic intervals composing that partition, where $1=i_1 <i_2<\ldots<i_{D_m}<i_{D_m+1}=n+1$. For all $1\leq k\leq D_m$, any column of $\hat{s}_m$ whose index belongs to $\{i_k,\ldots,i_{k+1}-1\}$ is equal to the mean $\bar{X}(i_k:i_{k+1})$ of the columns of $X$ whose indices belong to the interval $\{i_k,\ldots,i_{k+1}-1\}$. Owing to the form of the penalty, and to the additivity of the least-squares criterion, the whole criterion to minimize breaks down into a sum:
\begin{equation}
\|X-\hat{s}_m\|^2+\pen(m)=\sum_{k=1}^{D_m} \mathcal{L}(i_k,i_{k+1}),
\label{eq:additivecriterion}
\end{equation}
where, for all $1\leq k\leq D_m$,
$$\mathcal{L}(i_k,i_{k+1})=\cdi
+\sum_{i=i_k}^{i_{k+1}-1}\|X_i-\bar{X}(i_k:i_{k+1})\|^2_r.$$
By comparison with the method suggested in the previous paragraph, we are left with only one minimization problem, with no dimension constraint, instead of $n$.

We now turn to graph theory where our minimization problem finds a natural interpretation. We consider the weighted directed graph $G$ having $\{1,\ldots,n+1\}$ as vertex set and whose edges are the pairs $(i,j)$ such that $\{i,\ldots,j-1\}$ is a dyadic interval of $\{1,\ldots,n\}$ assigned with the weight $\mathcal{L}(i,j)$. A little vocabulary will be helpful. We say that a vertex $j$ is a successor to a vertex $i$ if $(i,j)$ is an edge of the graph $G$ and we associate to each vertex $i$ its successor list $\Gamma_i$. For all $1\leq D\leq n$, a $D+1$-uple $(i_1,i_2,\ldots,i_{D+1})$ of vertices of $G$ such that $i_1=1$, $i_{D+1}=n+1$ and each vertex is a successor to the previous one, will be called a path leading from $1$ to $n+1$ in $D$ steps. The length of such a path is defined as $\sum_{k=1}^{D} \mathcal{L}(i_k,i_{k+1}).$
Determining $\hat{m}$ thus amounts to finding a shortest path leading from $1$ to $n+1$ in the graph $G$. That problem can be solved by using one of the simplest shortest-path algorithms, the one dedicated to acyclic directed graphs, presented in~\cite{Cormenetal} (Section 24.2) for instance. For the sake of completeness, we also describe it in Table~\ref{flottants:algorithm:figure}. We have to underline that there are only $2n-1$ dyadic intervals of $\{1,\ldots,n\}$. Therefore, the graph $G$, with $n+1$ vertices and $2n-1$ edges, can be represented by only $\mathcal{O}(n)$ data: the weights $\mathcal{L}(i,j)$, for $1\leq i\leq n$ and $ j\in\Gamma_i$, and the successor lists $\Gamma_i$, for $1\leq i\leq n$.
\begin{table}[h]
\caption{Algorithm for computing $\tilde{s}$}
\begin{minipage}{12 cm}
\noindent
\rule{\textwidth}{0.5pt}

\emph{\textbf{Step 1} : Initialization}

Set $d(1)=0$ and $p(1)=+\infty$.

For $i=2,\ldots,n+1$,

\qquad set $d(i)=+\infty$ and $p(i)=+\infty$.

\bigskip
\emph{\textbf{Step 2} : Determining the lengths of the shortest paths with origin 1}

For $i=1,\ldots,n$,

\qquad for $j \in \Gamma_i$,

\qquad\qquad if $d(j)> d(i)+ \mathcal{L}(i,j)$,

\qquad\qquad\qquad then do $d(j)\gets d(i)+\mathcal{L}(i,j)$ and $p(j)\gets i$.

\bigskip
\emph{\textbf{Step 3} : Determining a shortest path $P$ from 1 to $n+1$}

Set $pred = p(n+1)$ and $P=(n+1)$.

While $pred \neq +\infty$,

\qquad replace $P$ with the concatenation of $pred$ followed by $P$,

\qquad do $pred \gets p(pred)$.

\bigskip
\emph{\textbf{Step 4} : Computing the preliminary estimator}

Set $\tilde{D}=length(P)-1$.

For $k=1,\ldots,\tilde{D},$

\qquad for $i=P(k),\ldots,P(k+1)-1$,

\qquad\qquad set $\tilde{s}_i=\bar{X}(P(k):P(k+1))$.

\noindent
\rule{\textwidth}{0.5pt}
\end{minipage}
\label{flottants:algorithm:figure}
\end{table}
In the key step of the algorithm, i.e. step 2, each edge is only considered once. When the time comes to consider the edges with origin $i$, the variables $d(i)$ and $p(i)$ respectively contain the length of a shortest path from 1 to $i$ and a predecessor of $i$ in such a path. Just before the edge $(i,j)$, where $j\in\Gamma_i$, be processed, the variables $d(j)$ and $p(j)$ contain respectively the length of a shortest path leading from 1 to $j$ and a predecessor of $j$ in such a path, based solely on the edges that have already been encountered. Then dealing with the edge $(i,j)$ consists in testing whether the length of the path leading from 1 to $j$ can be shortened by going via $i$ and updating, if necessary, $d(j)$ and $p(j)$. What clearly appears from the above description of the algorithm is that its complexity is only \emph{linear} in the size $n$ of the data.

\section{Hybrid estimator}\label{sec:hybrid}
\begin{figure}[h]
    \begin{center}
        \includegraphics[width=7cm,height=5cm]{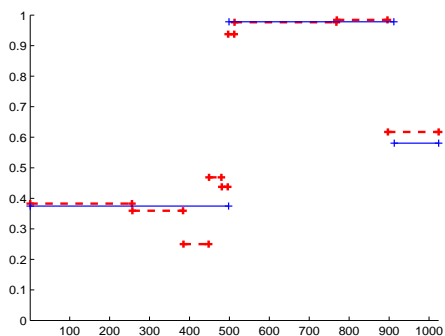}
    \caption{First lines of the distribution $s_a$ (full line) and of its preliminary estimator $\widetilde{s_a}$ (dashed line), as functions of $i$, $1\leq i \leq1024$.}
    \end{center}
    \label{flottants:testawithtildeopt:figure}
\end{figure}

Let us give a first glimpse of what can be expected from the preliminary estimator for detecting change-points in the distribution $s$. In figure~\ref{flottants:testawithtildeopt:figure}, we plot the first line of a distribution $s_a\in\mathscr{M}(2,1024)$ that is piecewise constant over a partition with only 3 segments together with the first line of a realization of $\widetilde{s_a}$. The value of \cdi\: has been chosen so as to minimize the distance between $s_a$ and its estimator.
\noindent
Both change-points in $s_a$ are indeed detected. But this example also shows that the selected partition, due to its special nature, is highly likely to contain some segments whose endpoints do not correspond to any significant rupture in $s$. In order to get rid of those, we propose a two-stage procedure, that we name hybrid procedure. After describing it, we provide an adaptivity result for that procedure and end this section with computational issues.

In the sequel, we suppose that $n \geq 2$.
In order to implement the hybrid procedure, we need to work with the set $\mathscr{M}(r,n/2)$ of $r\times (n/2)$ real matrices. That requires to define a series of notations, very close indeed to those encountered up to now. For all $t \in \mathscr{M}(r,n)$, we denote by $t^\bullet$ (resp. $t^\circ$) the element of $\mathscr{M}(r,n/2)$ composed of the columns of $t$ whose indices are even (resp. odd). We equip $\mathscr{M}(r,n/2)$ with the norm analogous to the norm $\|.\|$ on $\mathscr{M}(r,n)$. For the sake of simplicity, we will also denote by $\|.\|$ that norm on $\mathscr{M}(r,n/2)$. For a partition $m$ of $\{1,\ldots,n/2\}$, we denote by $S'_m$ the linear subspace of $\mathbb{R}^{n/2}$ generated by the indicator functions of the intervals $I\in m$ and by $D'_m$ its dimension.
These notations being settled, we are now able to define the hybrid estimator of $s$.
First, we compute the preliminary estimator of $s^\bullet$ based on $X^\bullet$, that is $\widetilde{s^\bullet}$, and we thus get a random partition of $\{1,\ldots,n/2\}$ into dyadic intervals denoted by $\hat{m}^\bullet$. Then, we consider the random collection $\widehat{\mathcal{M}}^\bullet$ of all the partitions of $\{1,\ldots,n/2\}$ that are built on $\hat{m}^\bullet$. For each partition $m$ of $\{1,\ldots,n/2\}$ into intervals, the least-squares estimators of $s^\circ$ in $\mathbb{R}^r \otimes S'_m$ is defined by
\begin{equation}
\widehat{s^\circ}_m = \argmin{t \in \mathbb{R}^r \otimes S'_m}\|X^\circ-t\|^2.
\notag
\end{equation}
We select
$$\hat{m}^\circ=\argmin{m \in \widehat{\mathcal{M}}^\bullet} \big\{\|X^\circ-{\widehat{s^\circ}}_m\|^2+\widehat{\pen}^{\circ}(m)\big\},$$
where the penalty $\widehat{\pen}^{\circ}$ will be chosen in the next paragraph. We define the penalized estimator of $s^\circ$ based on the collection $\widehat{\mathcal{M}}^\bullet$ as $\widehat{s^\circ}_{\hat{m}^\circ}.$
Last, we define the hybrid estimator $\tilde{s}_{hyb}$ of $s$ as the random matrix in $\mathscr{M}(r,n)$ whose submatrices composed respectively of columns with even indices and of columns with odd indices are both equal to $\widehat{s^\circ}_{\hat{m}^\circ}.$

Let us study $\tilde{s}_{hyb}$ from a theoretical point of view. Under a mild assumption on $s$, we derive from the results proved in the previous section the following adaptivity property for $\tilde{s}_{hyb}$.
\begin{theo}\label{Hybrid}
Let $\widehat{D}$ be the cardinal of $\hat{m}^\bullet$ and $\widehat{\pen}^{\circ} : \widehat{\mathcal{M}}^\bullet\rightarrow \mathbb{R}^+$ be a penalty of the form
\begin{equation}
\widehat{\pen}^{\circ}(m)=\Big(c_1 + c_2\ln\big(\widehat{D}/D'_m\big)\Big)D'_m,
\label{eq:penhyb}
\end{equation}
where $c_1$ and $c_2$ are positive. If $\|s^\bullet-s^\circ\|^2 \leq C\ln(n)$ and \cdi\:, $c_1$ and $c_2$ are large enough, then, for all $p\in (0,2],\alpha >1/p-1/2$ and $R$ such that $n^{-1/2}\leq R < n^\alpha$,
\begin{equation}
\sup_{s \in \mathscr{P}(\alpha, p, R)}\mathbb{E}_{s} \big[\|s - \tilde{s}_{hyb}\|^2\big]
\leq C\ln(n)\mathcal{R}(\alpha,p,R),
\label{eq:AdaptHybrid}
\end{equation}
where C only depends on $\cdi,c_1,c_2,\alpha$ and p.
\end{theo}
\noindent
Thus, with Inequality~\eqref{eq:AdaptHybrid}, we recover a result similar to Inequality~\eqref{eq:AdaptDyad}, up to a logarithmic factor.
\begin{proof}
For all $1 \leq D \leq \widehat{D}$, the number $\widehat{N}_D$ of partitions in $\widehat{\mathcal{M}}^\bullet$ with $D$ pieces satisfies
$$\widehat{N}_D= \binom{\widehat{D}-1}{D-1} \leq \bigg(\frac{e\widehat{D}}{D}\bigg)^D.$$
The above inequality results from a property of binomial coefficients that may be found in~\cite{Massart} (Proposition 2.5) for instance. So the weights defined by
\begin{equation}
L(D)=\ln (2e)+\ln (\widehat{D}/D), \text { for } 1 \leq D \leq \widehat{D},
\notag
\end{equation}
are such that
$$\sum_{D=1}^{\widehat{D}}\widehat{N}_D\exp(-DL(D)) \leq 1.$$
Moreover, the penalty $\widehat{\pen}^{\circ}$ given by~\eqref{eq:penhyb} fulfills the hypotheses of Theorem 1 in~\cite{DLT} provided $c_1$ and $c_2$ are large enough. With a slight abuse of notation, for any partition $m$ of $\{1,\ldots,n/2\}$, we still denote by $\bar{t}_m$ the orthogonal projection of an element $t\in\mathscr{M}(r,n/2)$ on $\mathbb{R}^r\otimes S'_m$.
Working conditionally to $X^\bullet$, the collection $\widehat{\mathcal{M}}^\bullet$ is deterministic, so we deduce from Theorem 1 of~\cite{DLT} applied to the estimator $\widehat{s^\circ}_{\hat{m}^\circ}$ of $s^\circ$ that
\begin{equation}
\mathbb{E}_{s^\circ} \big[\|s^\circ - \widehat{s^\circ}_{\hat{m}^\circ}\|^2|X^\bullet\big]
\leq C(c_1,c_2)\big[\|s^\circ - \overline{s^\circ}_{\hat{m}^\bullet}\|^2+\widehat{\pen}^\circ(\hat{m}^\bullet)\big].
\label{eq:riskcirc}
\end{equation}
We recall that $\widetilde{s^\bullet}=\widehat{s^\bullet}_{\hat{m}^\bullet}$. So, thanks to the triangle inequality, and since an orthogonal projection is a shrinking map, we get
$$\|s^\circ - \overline{s^\circ}_{\hat{m}^\bullet}\|^2
\leq C\big(\|s^\circ - s^\bullet\|^2 + \|s^\bullet - \widetilde{s^\bullet}\|^2\big).$$
Besides, for all $m \in \widehat{\mathcal{M}}^\bullet$,
$$\widehat{\pen}^\circ(m) \leq C(c_1,c_2)\ln(n)D'_m.$$
Taking into account the last two inequalities and integrating with respect to $X^\bullet$ then leads from~\eqref{eq:riskcirc} to
\begin{equation}
\mathbb{E}_{s} \big[\|s^\circ - \widehat{s^\circ}_{\hat{m}^\circ}\|^2\big]
\leq C(c_1,c_2)\Big[ \|s^\circ - s^\bullet\|^2 + \mathbb{E}_{s^\bullet} \big[\|s^\bullet - \widetilde{s^\bullet}\|^2\big]+\ln(n)\mathbb{E}_{s^\bullet}(D'_{\hat{m}^\bullet})\Big],
\notag
\end{equation}
where $D'_{\hat{m}^\bullet}$ is nothing but $\widehat{D}$.
Besides, it follows from the definition of $\tilde{s}_{hyb}$ that
\begin{equation*}
\|s-\tilde{s}_{hyb}\|^2
=\|s^\bullet-\widehat{s^\circ}_{\hat{m}^\circ}\|^2
+\|s^\circ-\widehat{s^\circ}_{\hat{m}^\circ}\|^2.
\end{equation*}
Applying the triangle inequality, we then get
\begin{equation}
\|s-\tilde{s}_{hyb}\|^2 \leq C\big(\|s^\bullet-s^\circ\|^2+\|s^\circ-\widehat{s^\circ}_{\hat{m}^\circ}\|^2\big).
\notag
\end{equation}
Consequently,
\begin{equation}
\mathbb{E}_{s} \big[\|s - \tilde{s}_{hyb}\|^2\big]
\leq C(c_1,c_2)\Big[ \|s^\circ - s^\bullet\|^2 + \mathbb{E}_{s^\bullet} \big[\|s^\bullet - \widetilde{s^\bullet}\|^2\big]+\ln(n)\mathbb{E}_{s^\bullet}(\widehat{D})\Big].
\label{eq:hybbtoodd}
\end{equation}
Let us denote by $\mathcal{M}'$ the set of all partitions of $\{1,\ldots,n/2\}$ into dyadic intervals. For the risk of $\widetilde{s^\bullet}$, Theorem 1 of~\cite{DLT} provides
\begin{equation}
\mathbb{E}_{s^\bullet} \big[\|s^\bullet- \widetilde{s^\bullet}\|^2\big]
\leq C(\cdi)\inf_{m \in \mathcal{M}'}\big\{\|s^\bullet - {\overline{s^\bullet}}_{m}\|^2+ D'_{m}\big\}.
\label{eq:riskspoint}
\end{equation}
In order to bound the term $\mathbb{E}_{s^\bullet}(\widehat{D})$, we need to go back to the proof of Theorem 1 in~\cite{DLT} (Section 8.1). As already seen during the proof of Proposition~\ref{diOracle}, we can choose a positive constant $L$ such that $\sum_{m\in\mathcal{M}'}\exp(-LD'_m) \leq 1$. Let us fix a partition $m\in \mathcal{M}'$ and $\xi > 0$. Using the same notation as in~\cite{DLT}, we deduce from the proof of Theorem 1 in~\cite{DLT} that there exists an event $\Omega_\xi(m)$ such that $\mathbb{P}_{s^\bullet}(\Omega_\xi(m)) \geq 1-\exp(-\xi)$ and on which
\begin{equation}
\cdi \widehat{D}
\leq C_1\|s^\bullet-\overline{s^\bullet}_m\|^2+C_2(\cdi)D'_m+C_3 \widehat{D}+C_4\xi.
\notag
\end{equation}
Therefore, if $\cdi > C_3$, then
\begin{equation}
\widehat{D}
\leq C(\cdi)\big(\|s^\bullet-\overline{s^\bullet}_m\|^2+D'_m+\xi\big).
\notag
\end{equation}
Integrating this inequality and taking the infimum over $m \in \mathcal{M}'$ then yields
\begin{equation}
\mathbb{E}_{s^\bullet}(\widehat{D})
\leq C(\cdi) \inf_{m \in \mathcal{M}'}\big\{\|s^\bullet-\overline{s^\bullet}_m\|^2+D'_m\big\}.
\label{eq:expectpenbullet}
\end{equation}
Moreover, one can check that
\begin{equation}
\inf_{m \in \mathcal{M}'}\big\{\|s^\bullet - {\overline{s^\bullet}}_{m}\|^2+ D'_{m}\big\}
\leq \inf_{m \in \mathcal{M}}\big\{\|s - \bar{s}_{m}\|^2+ D_{m}\big\}.
\label{eq:fromM'toM}
\end{equation}
Combining Inequalities~\eqref{eq:hybbtoodd} to~\eqref{eq:fromM'toM} and the assumption on $\|s^\bullet-s^\circ\|^2$, we finally get
\begin{equation}
\mathbb{E}_{s} \big[\|s- \tilde{s}_{hyb}\|^2\big]
\leq C(\cdi,c_1,c_2) \ln(n)\inf_{m \in \mathcal{M}}\big\{\|s - \bar{s}_m\|^2+ D_m\big\}.
\label{eq:riskhyb}
\notag
\end{equation}
We then conclude the proof as that of Theorem~\ref{adaptivity}.
\end{proof}

Regarding the computation of $\tilde{s}_{hyb}$, we know from Section~\ref{sec:computing} that determining $\widetilde{s^\bullet}$ only requires $\mathcal{O}(n)$ computations. On the other hand, since $\widehat{\pen}^\circ$ is not linear in the dimension of the models, $\hat{m}^\circ$ has to be determined following the method suggested at the beginning of Section~\ref{sec:computing} and using Bellman's algorithm. If we impose an upper-bound $D_{max}$ on the dimension of the model selected during the second stage, determining $\hat{m}^\circ$ given $X^\bullet$ then requires of the order of $\widehat{D}^2D_{max}$ computations. Since $\widehat{D}$ is upper-bounded by $n/2$, we can only ensure that the computational complexity of $\tilde{s}_{hyb}$ is, in the worst case, of the order of $n^2D_{max}$. However, we will see in Section~\ref{sec:simulations} that, in practice, the hybrid procedure can also be implemented with a linear complexity only and with quite satisfactory results.

\section{Simulation study}\label{sec:simulations}

In the previous sections, we were only interested in giving a form of penalty yielding, in theory, a performant estimator. The aim of this section is to study practical choices of the penalty for each procedure. Several simulations allow to assess the relevance of these choices and to illustrate the qualities of each procedure.

                \subsection{Choosing the penalty constant for the preliminary estimator}\label{sec:dipenalty}

\begin{figure}[p]
    \begin{center}
        \includegraphics[width=12cm,height=10cm]{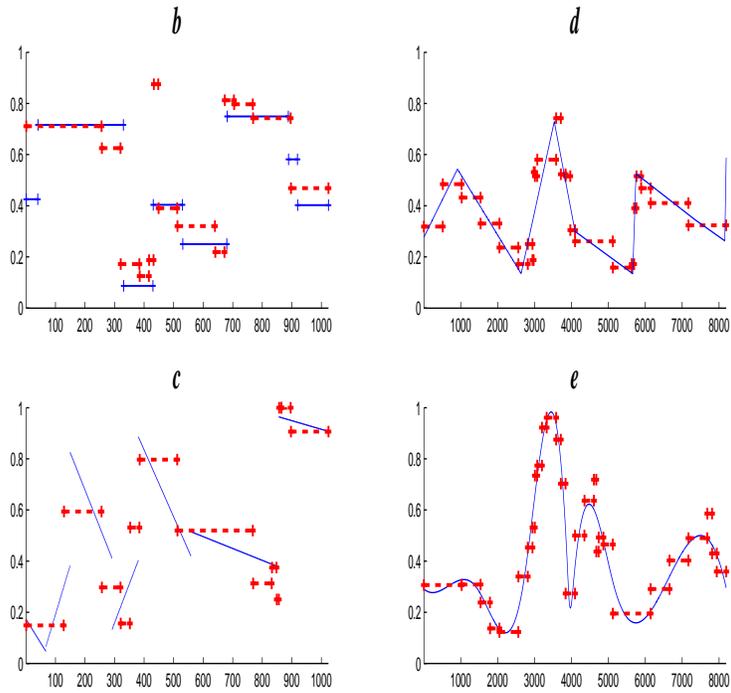}
    \end{center}
        \caption{\small{First lines of $s$ (full line) and $\tilde{s}$ (dashed line), computed with a data-driven penalty, for $s\in\{s_b,s_c,s_d,s_e\}$}.}
    \label{flottants:testbtoewithstilderand:figure}
\end{figure}

\begin{figure}[p]
    \begin{center}
        \includegraphics[width=12cm,height=6cm]{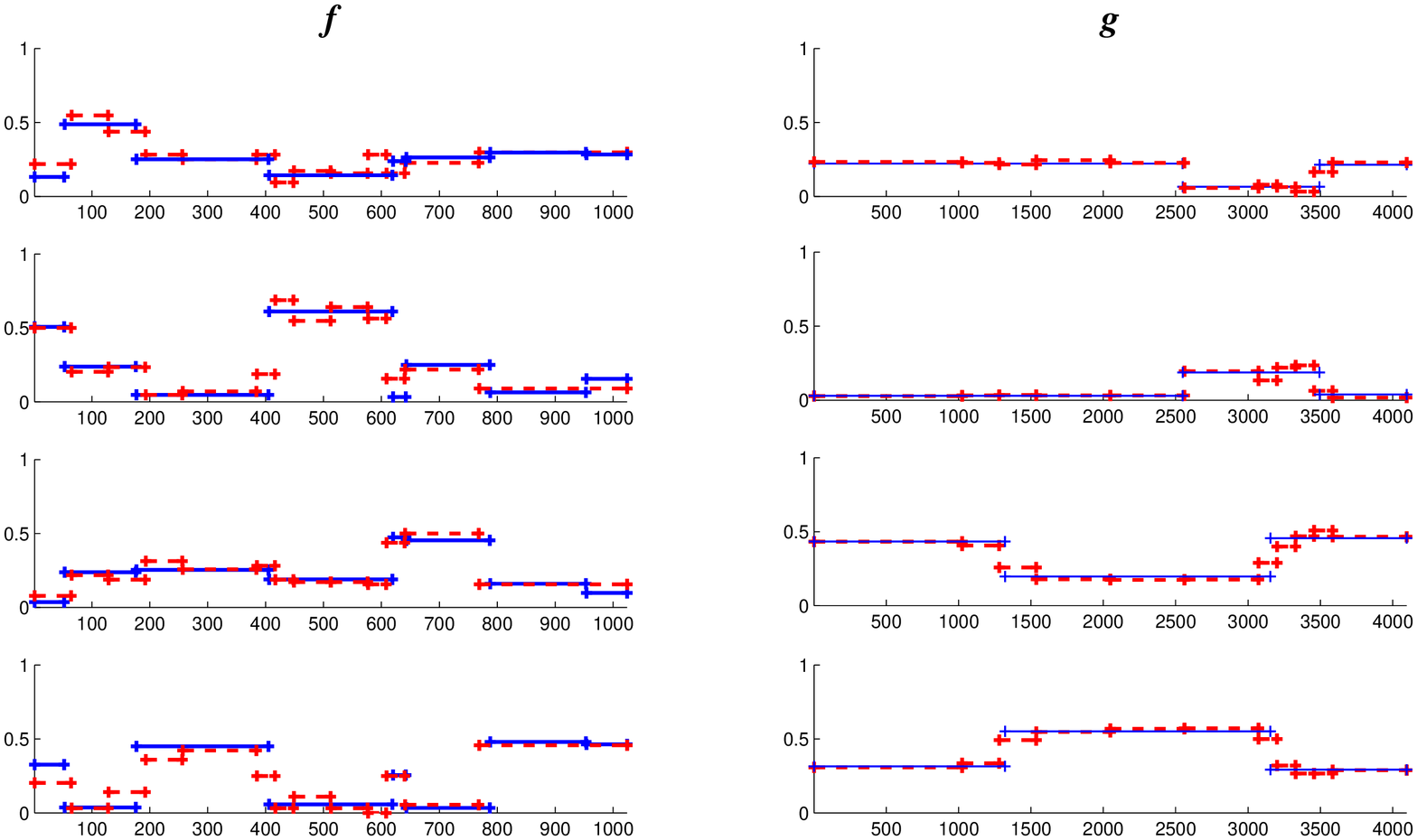}
        \bigskip
        \includegraphics[width=6cm,height=6cm]{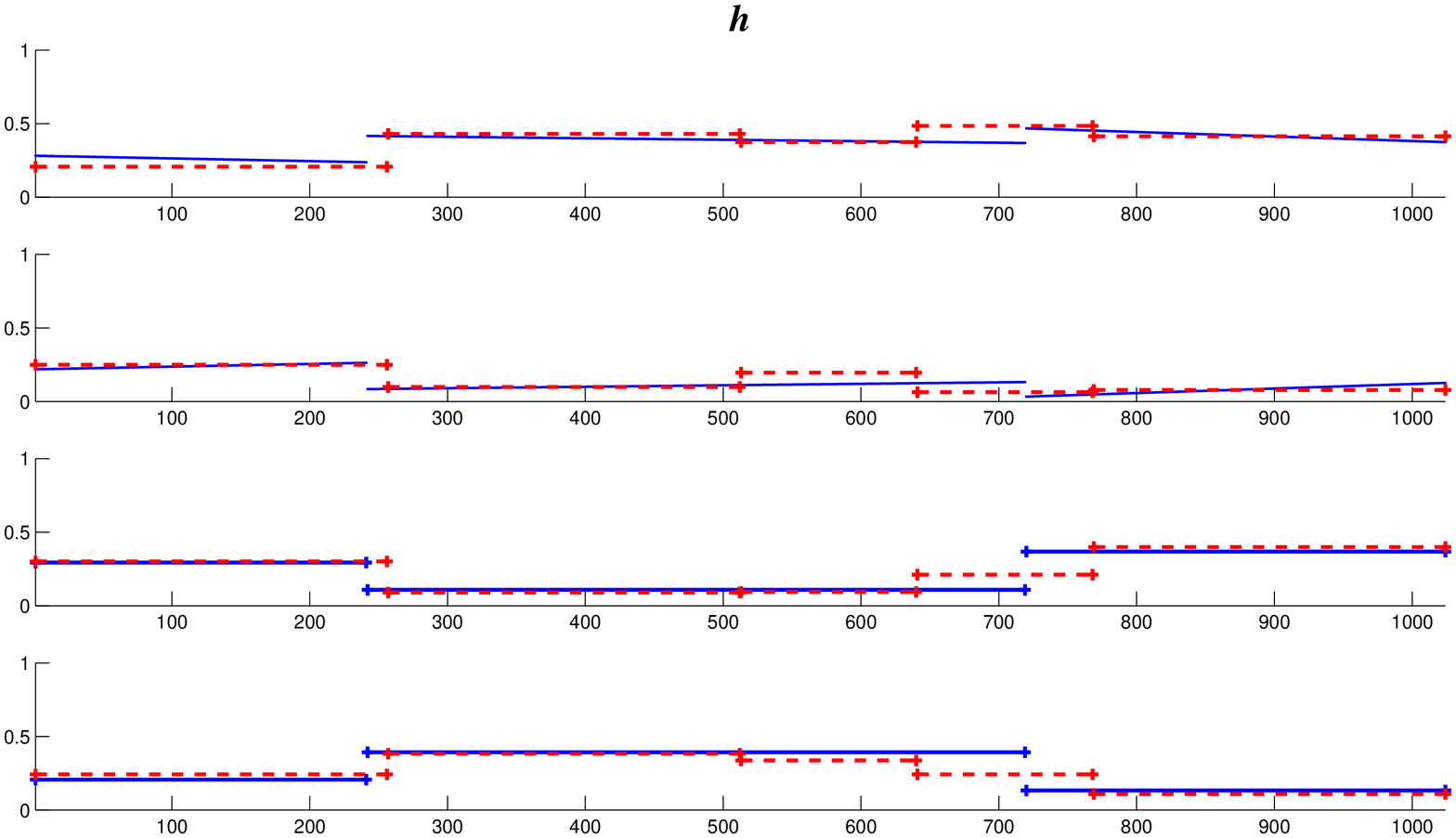}
    \end{center}
    \caption{\small{Four lines of $s$ (full line) and $\tilde{s}$ (dashed line), computed with a data-driven penalty, for $s\in \{s_f,s_g,s_h\}$}.}
    \label{flottants:testfandgwithstilderand:figure}
\end{figure}

We have examined the cases $r=2$ and $r=4$, with different values of $n=2^N$. For $r=2$, the distribution $s$ is entirely determined by its first line, that is the only one to be plotted, as a function of the parameter $i$, $1\leq i\leq n$ (cf. Figure~\ref{flottants:testawithtildeopt:figure} for $s_a$ and Figure~\ref{flottants:testbtoewithstilderand:figure} for $s_b$ to $s_e$). For $r=4$, examples $s_f$ to $s_h$ are plotted in Figure~\ref{flottants:testfandgwithstilderand:figure}. Part of our examples, $s_a$, $s_b$, $s_f$ and $s_g$, are piecewise constant. We also extend our study to other examples of distributions having jumps, such as $s_c$ and $s_h$, whose lines are piecewise affine. But the estimation capacities of $\tilde{s}$, and not only its ability to detect change-points, deserve to be illustrated. So, we also present smoother examples, if we may say so for functions of a discrete parameter, such as $s_d$ or $s_e$.

As already said in Section~\ref{sec:def}, the estimator $\tilde{s}$ has been designed for satisfying an oracle inequality, what it almost does according to Proposition~\ref{diOracle}. Therefore, the risk of the oracle, i.e. $\inf_{m \in \mathcal{M}}\mathbb{E}_s \big[\|s - \hat{s}_m\|^2\big]$, serves as a benchmark in order to judge of the quality of $\tilde{s}$, and also of the quality of a method for choosing a penalty constant.
We have studied two methods for choosing an adequate penalty constant. The different quantities introduced in the sequel have been estimated over 500 simulations.
The first method aims at determining the value of the constant \cdi\:that almost minimizes the risk of $\tilde{s}$, whatever $s$. Denoting by $\tilde{s}(c)$ the preliminary estimator when \cdi\:takes the value $c$, we have estimated
$$c^\star(s):=\argmin{c}\:\mathbb{E}_s \big[\|s - \tilde{s}(c)\|^2\big],$$
where, in practice, we have varied $c$ from 0 to 4, by step 0.1, and from 4 to 6 by step 0.5.
We plot in Table~\ref{flottants:testdi:table} an estimation of $c^\star$ and the ratio $Q^\star$ between an estimation of $\mathbb{E}_s \big[\|s - \tilde{s}(c^\star)\|^2\big]$ and the estimated risk of the oracle.
In view of the results obtained here, we come to the following conclusions: taking $\cdi=2$ seems reasonable when $r=2$, but taking $\cdi=2.5$ seems more appropriate when $r=4$.
We give in Table~\ref{flottants:testdi:table} the ratio $Q_{c}$ between the estimated risk of $\tilde{s}(c)$ and the estimated risk of the oracle, where $c=2$ for $r=2$ and $c=2.5$ for $r=4$. Comparing $Q_{c}$ to $Q^\star$ confirms that the choice of those values for is \cdi\:relevant.
Nevertheless, a good penalty should adapt to the unknown distribution $s$ to estimate. That's why we have also tried a data-driven method, inspired from results proved by Birgé and Massart in a Gaussian framework (cf.~\cite{BirgéMassart2}). That method has already been implemented in the same framework as ours in~\cite{DLT}, Section 8. Given a simulation of $(Y_1,\ldots,Y_n)$, the procedure we have followed can be decomposed in three steps:
\begin{itemize}
\item [\textbullet] determine the dimension $\hat{D}(c)$ of the selected partition for each value $c$ of the penalty constant \cdi, where one varies $c$ from 0 to 3, by step 0.1;
\item [\textbullet] compute the difference between the dimensions of the selected partitions for two consecutive values of \cdi\: and retain the value $\hat{c}$ corresponding to the biggest jump in dimension under the constraint $\hat{D}(\hat{c}) \leq D_{max}$, where $D_{max}$ is a prescribed  maximal dimension;
\item [\textbullet] choose the constant $\hat{c}_{j}=2\hat{c}$ to compute the preliminary estimator.
\end{itemize}
We have taken $D_{max}$ of the order of $n/(\ln(n))^\mu$, with $\mu$ close to 2. That choice is inspired in fact both from the method proposed in~\cite{SzpanSzpanRen} and from a constraint appearing in the theoretical results of~\cite{LebarbierNédélec} when using a penalized maximum likelihood criterion (cf. Condition (2.17) in Theorem 2.3. of~\cite{LebarbierNédélec}). That choice seems to yield good results, whatever $s$ or $n$. Here we have set $D_{max}=30$ when $N=10$, $D_{max}=100$ when $N=12$ and $D_{max}=175$ when $N=13$. In order to assess the performance of that second method, we give in Table~\ref{flottants:testdi:table} the ratio $Q_{j}$ between the estimated risk of $\tilde{s}$ for that procedure and the estimated risk of the oracle. We also give estimations of the mean value and the standard-error of $\hat{c}_{j}$, denoted respectively by $\bar{c}_{j}$ and $\sigma_{j}$.

\begin{table}[h]
\caption{Performance of the preliminary estimator for different choices of the penalty constant.}
    \begin{center}
        \begin{tabular}{c||c|c||c|c|c||c|c|c}
$s$ &$r$ &$N$ &$c^\star$ &$Q^\star$ &$Q_{c}$ &$\bar{c}_{j}$ &$\sigma_{j}$ &$Q_{j}$\\ \hline
$s_a$ &$2$ &$10$ &$1.7$ &$2.4$ &$2.6$  &$2.2$ &$0.3$ &$2.7$\\
$s_b$ &$2$ &$10$ &$1.7$ &$1.9$ &$1.9$ &$2.6$ &$0.4$ &$2.1$\\
$s_c$ &$2$ &$10$ &$1.8$ &$1.8$ &$1.8$ &$2.5$ &$0.4$ &$1.8$\\
$s_d$ &$2$ &$13$ &$2.2$ &$1.5$ &$1.6$ &$2.2$ &$0.1$ &$1.6$\\
$s_e$ &$2$ &$13$ &$2.2$ &$1.7$ &$1.8$  &$2.2$ &$0.1$ &$1.8$ \\ \hline
$s_f$ &$4$ &$10$ &$2$ &$1.4$ &$1.5$  &$3.3$ &$0.6$ &$1.7$\\
$s_g$ &$4$ &$12$ &$2.5$ &$1.3$ &$1.4$  &$2.5$ &$0.1$ &$1.3$\\
$s_h$ &$4$ &$10$ &$2.6$ &$1.3$ &$1.3$  &$2.7$ &$0.2$ &$1.3$\\
        \end{tabular}
    \end{center}
\label{flottants:testdi:table}
\end{table}

Let us analyze the results of the simulations. In terms of risk, both methods have in fact roughly the same performance. Nevertheless, the first one requires to calibrate anew a constant when changing the value of $r$, whereas the data-driven method has the advantage to automically adapt to the value of $r$. Therefore, the latter should be recommended, and that is the one we have used to build the estimators plotted in Figures~\ref{flottants:testbtoewithstilderand:figure} and~\ref{flottants:testfandgwithstilderand:figure}. Let us now examine the values of $Q^\star$ (or $Q_{c}$, or $Q_{j}$) for the different examples. As foreseen by the oracle-type inequality~\eqref{eq:oracletypedi}, the ratio between the risk of the preliminary estimator and that of the oracle depends on $s$. In particular, the ratios $Q^\star$, $Q_{c}$ or $Q_{j}$ reach their highest value for $s_a$. It should be noted that the first line of this example takes values very close to 1 on a large segment (cf. Figure~\ref{flottants:testawithtildeopt:figure}), a critical case according to the oracle-type inequality. However, for all examples studied here, the values of those ratios remain quite low, inferior or close to 2, except for $s_a$.

                \subsection{Choosing the penalty constants for the hybrid estimator}\label{sec:hybpenalty}
For the first stage of the hybrid procedure, the preliminary estimator has been computed using the data-driven penalty.
For the second stage, the practical choice of an adequate penalty is more delicate, since the theoretical penalty depends in this case on two constants and on the dimension $\widehat{D}$ of the partition selected during the first stage. We have first tried here the same method as Lebarbier in~\cite{Lebarbier}, Chapter 7, for her own hybrid procedure. So we have assigned to all partitions of $\{1,\ldots,n/2\}$ into $D$ intervals the same penalty
$$\widehat{\pen}^{\circ}_1(D)=\hat{\beta}_1(2.5+\ln(\widehat{D}/D))D,$$
where $\hat{\beta}_1$ is determined according to the same process as $\hat{c}_{j}$.
That penalty is proportional to the penalty calibrated by Lebarbier in~\cite{Lebarbier} (Chapter 3). The latter was in fact designed for the estimation of a regression function in a Gaussian framework via model selection based on an exhaustive collection of partitions. Anyway, the major drawback of such a method, as said at the end of Section~\ref{sec:hybrid}, is that we are only able to evaluate its worst case computational complexity, of the order of $\mathcal{O}(n^3)$. So we have also tried to assign to all partitions of $\{1,\ldots,n/2\}$ into $D$ intervals the penalty
$$\widehat{\pen}^{\circ}_2(D)=\hat{\beta}_2D,$$
where $\hat{\beta}_2$ is determined once again according to the same process as $\hat{c}_{j}$. Since that penalty is a linear funtion of $D$, the hybrid procedure can be implemented in that case with only $\mathcal{O}(n)$ computations.

\begin{figure}[p]
    \begin{center}
        \includegraphics[width=12cm,height=10cm]{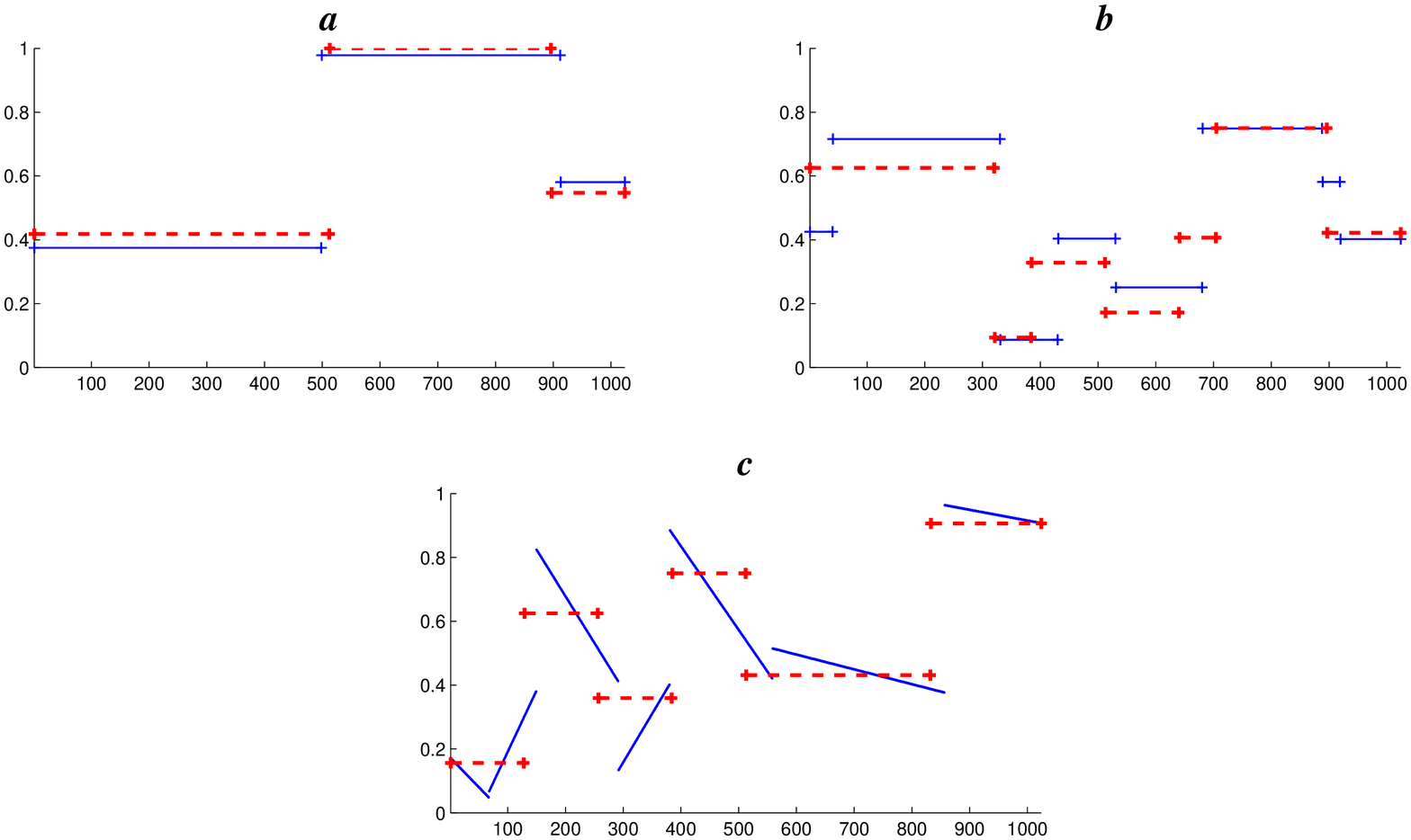}
    \end{center}
    \caption{\small{First lines of $s$ (full line) and of its hybrid estimator (dashed line) for $s \in \{s_a,s_b,s_c\}$}.}
    \label{flottants:Hybrid2:figure}
\end{figure}

\begin{figure}[p]
    \begin{center}
        \includegraphics[width=12cm,height=6cm]{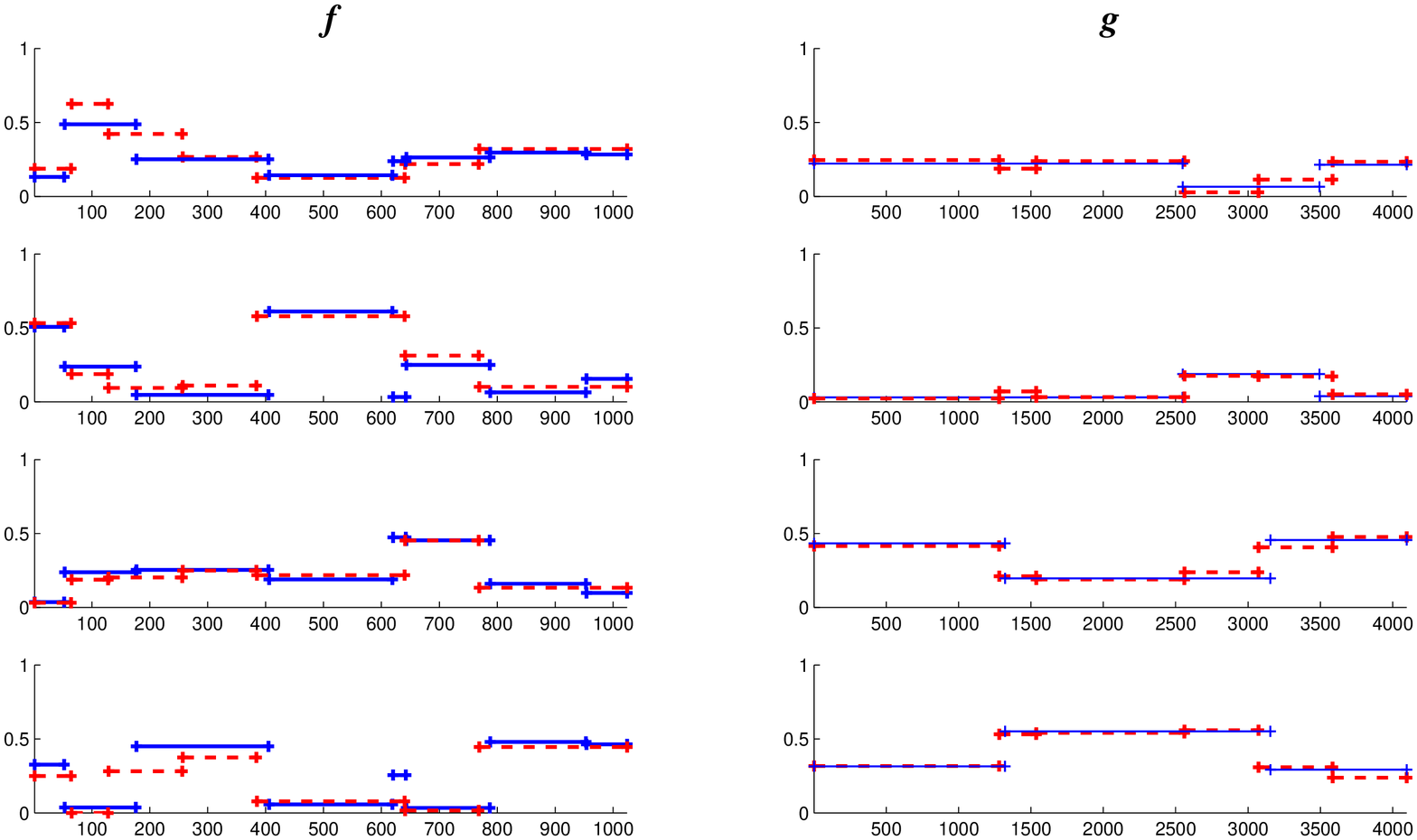}
        \bigskip
        \includegraphics[width=6cm,height=6cm]{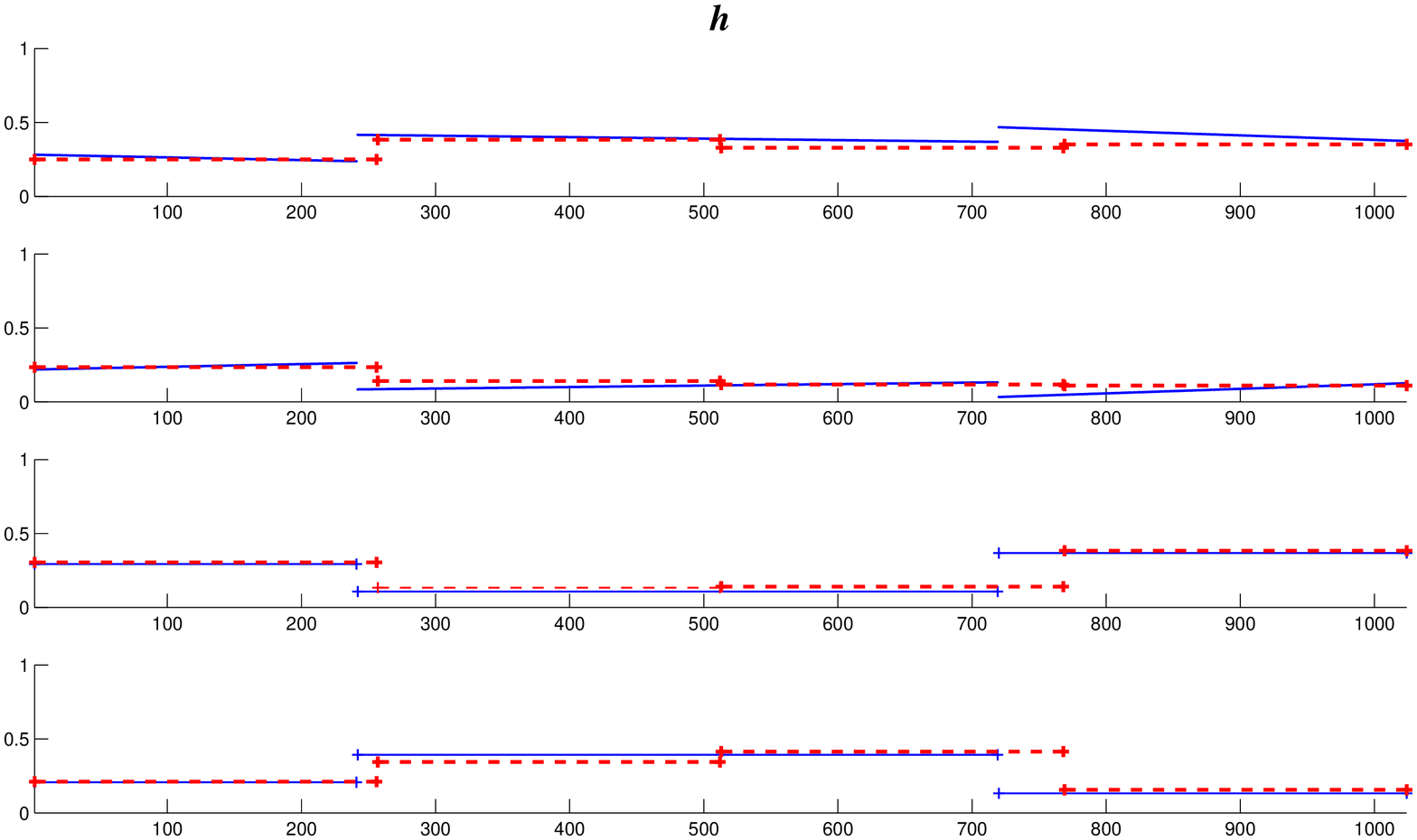}
    \end{center}
    \caption{\small{Four lines of $s_f$, $s_g$ and $s_h$ (full line) together with their hybrid estimators (dashed line).}}
    \label{flottants:Hybrid4:figure}
\end{figure}

In order to draw a comparison between these procedures and with the preliminary one, we give in Table~\ref{flottants:testhybrid:table} the following information for the distributions $s_a$ to $s_c$ and $s_f$ to $s_g$, still computed over 500 simulations. We first recall the dimension $D$ of the partition on which $s$ is built. Then we indicate the average dimensions $\bar{D}_0$ and $\bar{D}_i$ of the partitions selected respectively by the preliminary procedure, with a data-driven penalty, and the hybrid procedure with $\widehat{\pen}^{\circ}_i$, for $i\in\{1,2\}$. We also give the average value $Q_{i:0}$ of the ratio between the estimated risk of the hybrid estimator for $\widehat{\pen}^{\circ}_i$, for $i\in\{1,2\}$, and the estimated risk of the preliminary estimator. Let us compare both ways to implement the hybrid procedure. We observe that $Q_{2:0}$ is almost always of the same order as $Q_{1:0}$, and even slightly lower in most cases. Therefore, taking into account the computational complexity, we cannot but recommend to use $\widehat{\pen}^{\circ}_2$. That is the choice we have made for the hybrid estimators represented in Figures~\ref{flottants:Hybrid2:figure} and~\ref{flottants:Hybrid4:figure}. Let us now compare the hybrid procedure with the preliminary one for the examples under study. First, the values of $\bar{D}_2$ and $\bar{D}_0$ indicate that, with the former, the dimension of the selected partition is much closer to the true one. Moreover, the figures show that the most significant ruptures are still detected, are quite close to the true ones, and that irrelevant ruptures are much fewer with the hybrid procedure. The only price to pay is an increase in risk, but only by a factor of the order of 1.5.

\begin{table}[p]
\caption{Comparison between the hybrid procedure, for different penalties, and the preliminary procedure.}
    \begin{center}
        \begin{tabular}{c||c|c|c|c||c|c}
$s$ &$D$ &$\bar{D}_0$ &$\bar{D}_1$ &$\bar{D}_2$ &$Q_{1:0}$ &$Q_{2:0}$\\ \hline
$s_a$ &$3$ &$7.7$ &$3.0$ &$3.4$ &$1.4$ &$1.3$\\
$s_b$ &$8$ &$13.4$ &$4.9$ &$6.9$ &$1.5$ &$1.4$\\
$s_c$ &$7$ &$11.7$ &$5.1$ &$5.1$ &$1.7$ &$1.8$\\
$s_f$ &$8$ &$11.5$ &$4.2$ &$5.9$ &$2.1$ &$1.6$\\
$s_g$ &$5$ &$11.5$ &$7.8$ &$6.9$ &$1.8$ &$1.4$\\
$s_h$ &$3$ &$5.3$ &$2.2$ &$3.1$ &$2.3$ &$1.5$\\
        \end{tabular}
    \end{center}
\label{flottants:testhybrid:table}
\end{table}

\section{Proof of the approximation result}\label{sec:proof}

In this section, we prove Theorem~\ref{treeapprox} following the same path as DeVore and Yu in~\cite{DeVoreYu}. We first describe the approximation algorithm on which that result relies. Then, we give the main lines of the proof and also demonstrate the key result, that is a direct consequence of the approximation algorithm. The proofs of more technical points are postponed to the next subsections.

            \subsection{Approximation algorithm}\label{sec:approxalgo}

Let us fix $p\in(0,2]$, $\alpha > 1/p - 1/2$, $R > 0$ and $D \in\{1,\ldots,n\}$. In order to prove Theorem~\ref{treeapprox}, we look for an upper bound for
$$\inf_{m \in \mathcal{M}_D}\|t - \bar{t}_m\|^2$$
uniformly over $t\in\mathscr{B}(\alpha, p, R)$.
An element $t\in\mathscr{M}(r,n)$ being fixed, the adaptive approximation algorithm presented by DeVore and Yu in~\cite{DeVoreYu} allows to generate partitions into dyadic intervals depending on $t$ such that the approximation error over each interval of the partitions is lower than a prescribed threshold.
An adequate choice of that threshold is expected to yield a partition, depending on $t$, that belongs to $\mathcal{M}_D$ and almost realizes the above infimum. In order to describe precisely the algorithm and the way to use it for our approximation problem, let us introduce some notations. Let $I$ be a dyadic interval of $\{1,\ldots,n\}$. The restriction of the norm $\|.\|$ to $I$ is denoted by $\|.\|_I$. Let $U$ be the linear subspace of $\mathbb{R}^n$ generated by the vector $(1\ldots 1)$, we denote by $\mathcal{E}_2(t,I)$ the error in approximating $t$ on $I$ by an element of $\mathbb{R}^r\otimes U$, i.e.
$$\mathcal{E}_2(t,I) = \inf_{c\in\mathbb{R}^r\otimes U}\|t-c\|_{I}.$$
Besides, both intervals obtained by dividing $I$ into two intervals of same length are called the children of $I$.
The algorithm proceeds as follows. We fix a threshold $\epsilon > 0$.
At the beginning, the set $\mathcal{I}^1(t,\epsilon)$ contains $I_{(0,0)}=\{1,\ldots,n\}$.
If $\mathcal{E}_2\big(t,I_{(0,0)}\big) \leq \epsilon$, then the algorithm stops. Else, $I_{(0,0)}$ is replaced in the partition $\mathcal{I}^1(t,\epsilon)$ with his children, hence a new partition $\mathcal{I}^2(t,\epsilon)$ of $\{1,\ldots,n\}$.
In the same way, the $k$-th step starts with a partition $\mathcal{I}^k(t,\epsilon)$ of $\{1,\ldots,n\}$ into $k$ dyadic intervals. If $\sup_{I\in\mathcal{I}^k(t,\epsilon)}\mathcal{E}_2(t,I) \leq \epsilon$, then the algorithm stops, else an interval $I$ such that $\mathcal{E}_2(t,I) > \epsilon $ is chosen in $\mathcal{I}^k(t,\epsilon)$
and replaced with his children, hence a new partition $\mathcal{I}^{k+1}(t,\epsilon)$ of $\{1,\ldots,n\}$ into $k+1$ dyadic intervals.
The algorithm finally stops, giving a partition $\mathcal{I}(t,\epsilon)$. Denoting by $S(t,\epsilon)$ the linear space composed of the functions that are piecewise constant on $\mathcal{I}(t,\epsilon)$, the approximation $A(t,\epsilon)$ of $t$ associated with this partition is defined as the orthogonal projection of $t$ on $\mathbb{R}^r\otimes S(t,\epsilon)$. So, the approximation error of $t$ by $A(t,\epsilon)$ satisfies
\begin{equation}
\|t-A(t,\epsilon)\|^2
=\sum_{I\in\mathcal{I}(t,\epsilon)}\big(\mathcal{E}_2(t,I)\big)^2
\leq |\mathcal{I}(t,\epsilon)|\epsilon^2.
\label{eq:majtree}
\notag
\end{equation}
For any $\epsilon >0$ such that the algorithm stops at the latest at step $D$, the approximation of $t$ that we get belongs to the collection $\{\mathbb{R}^r\otimes S_m\}_{m \in \mathcal{M}_D}$. Therefore
\begin{equation}
\inf_{m \in \mathcal{M}_D}{\|t - \bar{t}_m\|^2}
\leq |\mathcal{I}(t,\epsilon)|\epsilon^2.
\label{eq:premsapprox}
\notag
\end{equation}
Let us denote by $\mathcal{E}_D(t)$ the infimum of $|\mathcal{I}(t,\epsilon)|\epsilon^2$ taken over all $\epsilon > 0$ satisfying $|\mathcal{I}(t,\epsilon)| \leq D.$ This is in fact the quantity that we shall bound, as indicated in Theorem~\ref{algoapprox} below.
\begin{theo}\label{algoapprox}
Let $p\in (0,2]$, $\alpha > 1/p-1/2$ and $R>0$. For all $D \in \{1,\ldots,n\}$ and $t\in \mathscr{B}(\alpha, p, R)$,
\begin{equation}
\mathcal{E}_D(t)
\leq C(\alpha,p)nR^2D^{-2\alpha}.
\notag
\end{equation}
\end{theo}
\noindent
We then get Theorem~\ref{treeapprox} as a straightforward consequence of Theorem~\ref{algoapprox}.

        \subsection{Proof of Theorem~\ref{algoapprox}: the main lines}\label{sec:prooftheo}

Here are the notions and notations that we will need along the proof.
Let $p>0$, $\alpha > 0$ and $t \in \mathscr{M}(r,n)$. For every subset $I$ of $\{1,\ldots,n\}$, let
$$\mathcal{E}_p(t,I)=\inf_{v\in\mathbb{R}^r}\bigg(\sum_{k\in I}\|t_k-v\|_r^p\bigg)^{1/p}.$$
We define the vector $t^{\sharp,\alpha,p}$ in $\mathbb{R}^n$ whose coordinates are
$$t^{\sharp,\alpha,p}_i=\sup_{I \ni i}|I|^{-(\alpha+1/p)}\mathcal{E}_p(t,I),\text{ for } i=1,\ldots,n,$$
where the supremum is taken over all the dyadic intervals $I$ of $\{1,\ldots,n\}$ that contain $i$.
We denote by $\|.\|_{\ell_p}$ the (quasi-)norm defined on $\mathbb{R}^n$ by
$$\|u\|_{\ell_p} = \bigg(\sum_{i=1}^n|u_i| ^p\bigg)^{1/p}$$
(that is a norm only for $p\geq 1$) and by $\|.\|_{\ell_p,I}$ its restriction to a subset $I$ of $\{1,\ldots,n\}$. We define on $\mathbb{R}^n$ the discrete Hardy-Littlewood maximal function $M_p$ by
$$\big(M_p(u)\big)_i=\sup_{I \ni i}|I|^{-1/p}\|u\|_{\ell_p,I},\text{ for }
i=1,\ldots,n,$$
where the supremum is taken over all the dyadic intervals $I$ of $\{1,\ldots,n\}$ containing $i$. Last, we recall that every vector $u\in\mathbb{R}^n$ is identified with the function defined on $\{1,\ldots,n\}$ whose value in $i$ is $u_i$, for $1\leq i\leq n$, hence the meaning of notations such as $u\leq v$ or $u^q$, where $u\in\mathbb{R}^n$, $v\in\mathbb{R}^n$ and $q>0$.

The beginning of the proof directly results from the way the algorithm works out. A dimension $D$ being fixed, choosing $\epsilon>0$ as small as possible such that the algorithm generates a partition with at most $D$ intervals leads to a first comparison between the quantity $\mathcal{E}_D(t)$ and $D^{-2\alpha}$, without making use of any particular hypothesis on $t$.
\begin{prop}\label{proptree}
Let $\alpha > 0$ and $p(\alpha)=(\alpha+1/2)^{-1}$. For all $D\in\{1,\ldots,n\}$ and $t\in \mathscr{M}(r,n)$,
\begin{equation}
\mathcal{E}_D(t)
\leq C(\alpha)\|t^{\sharp,\alpha,2}\|_{\ell_{p(\alpha)}}^{2}D^{-2\alpha}.
\notag
\end{equation}
\end{prop}
\begin{proof}
If $t^{\sharp,\alpha,2}=0$, then, whatever $\epsilon >0$, $\mathcal{E}_2\big(t,I_{(0,0)}\big)\leq \epsilon$, so $\mathcal{E}_D(t)=0$, which completes the proof in that case.
Let us now suppose that $t^{\sharp,\alpha,2}$ is non-null, and let $\epsilon >0$. If $\mathcal{E}_2\big(t,I_{(0,0)}\big)\leq \epsilon$, then $|\mathcal{I}(t,\epsilon)|=1$. Else, let $I$ be a dyadic interval that belongs to $\mathcal{I}(t,\epsilon)$, then $I$ is a child of a dyadic interval $\tilde{I}$ such that
\begin{equation}
\epsilon < \mathcal{E}_2\big(t,\tilde{I}\big).
\notag
\end{equation}
Using the definition of $t^{\sharp,\alpha,2}$, we get, for all $i \in \tilde{I}$,
\begin{equation}
\mathcal{E}_2\big(t,\tilde{I}\big) \leq \big|\tilde{I}\big|^{\alpha+1/2}t^{\sharp,\alpha,2}_i.
\notag
\end{equation}
Since $I \subset\tilde{I}$, $|\tilde{I}|=2|I|$ and $p(\alpha)=(\alpha + 1/2)^{-1}$, the last two inequalities lead, for all $i \in I$, to
\begin{equation}
\epsilon < 2^{1/p(\alpha)}|I|^{1/p(\alpha)}t^{\sharp,\alpha,2}_i,
\notag
\end{equation}
hence
\begin{equation}
{\epsilon}^{p(\alpha)} < 2\sum_{i\in I}{\big(t^{\sharp,\alpha,2}_i\big)}^{p(\alpha)}.
\notag
\end{equation}
Then we deduce by summing over all the intervals $I$ in the partition $\mathcal{I}(t,\epsilon)$ that
\begin{equation}
|\mathcal{I}(t,\epsilon)|
\leq 2{\|t^{\sharp,\alpha,2}\|}_{\ell_{p(\alpha)}}^{p(\alpha)}{\epsilon}^{-p(\alpha)}.
\notag
\end{equation}
Whether $\mathcal{E}_2\big(t,I_{(0,0)}\big)\leq \epsilon$ or not, by choosing
$\epsilon=2^{1/p(\alpha)}{\|t^{\sharp,\alpha,2}\|}_{\ell_{p(\alpha)}}D^{-1/{p(\alpha)}}$, we get a partition $\mathcal{I}(t,\epsilon)$ that contains at most $D$ elements and satisfies
\begin{align}
|\mathcal{I}(t,\epsilon)|{\epsilon}^2
& \leq D^{1-2/p(\alpha)}2^{2/p(\alpha)}{\|t^{\sharp,\alpha,2}\|}_{\ell_{p(\alpha)}}^{2}.
\notag
\end{align}
As $p(\alpha)=(\alpha + 1/2)^{-1}$, we conclude that
\begin{align}
|\mathcal{I}(t,\epsilon)|{\epsilon}^2
& \leq 4^{\alpha+1/2}{\|t^{\sharp,\alpha,2}\|}_{\ell_{p(\alpha)}}^{2}D^{-2\alpha}.
\notag
\end{align}
\end{proof}

The proof of Theorem~\ref{algoapprox} now relies upon three inequalities. The first one allows to draw a comparison between $\mathcal{E}_D(t)$ and $D^{-2\alpha}$ via a term that does not depend on $t^{\sharp,\alpha,2}$ anymore but on $t^{\sharp,\alpha,p(\alpha)}$. It is the discrete analogue of a particular case of Theorem 4.3. of~\cite{DeVShar}.
\begin{prop}\label{changenorm}
Let $\alpha > 0$ and $p(\alpha)=(\alpha + 1/2)^{-1}$.
For all $t\in \mathscr{M}(r,n)$,
\begin{equation}
t^{\sharp,\alpha,2} \leq C(\alpha) M_{p(\alpha)}\big(t^{\sharp,\alpha,p(\alpha)}\big).
\notag
\end{equation}
\end{prop}
\noindent
From Propositions~\ref{proptree} and~\ref{changenorm}, we easily deduce that, for $\alpha > 0$, $p(\alpha)=(\alpha + 1/2)^{-1}$ and $D \in \{1,\ldots,n\}$,
\begin{equation}
\mathcal{E}_D(t)
\leq C(\alpha)\big\|M_{p(\alpha)}\big(t^{\sharp,\alpha,p(\alpha)}\big)\big\|_{\ell_{p(\alpha)}}^{2}D^{-2\alpha}.
\label{eq:prems}
\notag
\end{equation}
Let us now fix $p \in (0,2].$ By Jensen's inequality, we have
\begin{equation}
\big\|M_{p(\alpha)}\big(t^{\sharp,\alpha,p(\alpha)}\big)\big\|_{\ell_{p(\alpha)}}
\leq n^{1/p(\alpha)-1/p}\big\|M_{p(\alpha)}\big(t^{\sharp,\alpha,p(\alpha)}\big)\big\|_{\ell_{p}}
\notag
\end{equation}
and
\begin{equation}
t^{\sharp,\alpha,p(\alpha)}
\leq t^{\sharp,\alpha,p},
\notag
\end{equation}
hence
\begin{equation}
\mathcal{E}_D(t)
\leq C(\alpha) n^{2(\alpha+1/2-1/p)}
\big\|M_{p(\alpha)}(t^{\sharp,\alpha,p})\big\|_{\ell_{p}}^{2}D^{-2\alpha}.
\label{eq:lpnormtsharp}
\notag
\end{equation}
Though the most obvious comparison between a vector $u$ and any of its maximal functions is that the latter are greater than the first, the following maximal inequality also ensures a control of $u$ over its maximal functions (cf. inequality~\eqref{eq:genmaxineq} below). That inequality is in fact the discrete version of a fundamental result in functional analysis, namely the Hardy-Littlewood maximal inequality, that may be found in~\cite{BennShar} (Theorem 3.10) for instance.
\begin{prop}\label{maxf}
Let $q>1$. For all $u\in \mathbb{R}^n$,
\begin{equation}
\|M_1(u)\|_{\ell_q}
\leq C(q) \|u\|_{\ell_q}.
\notag
\end{equation}
\end{prop}
\noindent
Since the maximal function $M_q$, $q>0$, is related to $M_1$ by the property
\begin{equation}
M_q(u)= \big(M_1(u^q)\big)^{1/q}, \text{for all } u \in \mathbb{R}^n,
\notag
\end{equation}
Proposition~\ref{maxf} yields, for all $r>q>0$ and $u \in \mathbb{R}^n$,
\begin{equation}
\|M_q(u)\|_{\ell_r}\leq C(r,q) \|u\|_{\ell_r}.
\label{eq:genmaxineq}
\end{equation}
Thus, when applied with $u=t^{\sharp,\alpha,p}$, $r=p$ and $q=p(\alpha)$, this inequality leads to
\begin{equation}
\mathcal{E}_D(t)
\leq C(\alpha,p)n^{2(\alpha+1/2-1/p)}\|t^{\sharp,\alpha,p}\|_{\ell_{p}}^2D^{-2\alpha}.
\notag
\end{equation}
Last, Proposition~\ref{compSharpBesov} below provides the adequate control of the $\ell_p$-(quasi-)norm of $t^{\sharp,\alpha,p}$ by the size of the wavelet coefficients of $t$ and allows to complete immediately the proof of Theorem~\ref{algoapprox}.
\begin{prop}\label{compSharpBesov}
Let $p\in(0,2]$ and $\alpha > 1/p-1/2$.
For all $t \in \mathscr{M}(r,n)$,
$$\|t^{\sharp,\alpha,p}\|_{\ell_{p}}
\leq C(\alpha,p)n^{-(\alpha+1/2-1/p)}\Bigg(\sum_{j=0}^{N-1} 2^{jp(\alpha + 1/2 - 1/p)}\sum_{\lambda \in \Lambda(j)}\|\beta_\lambda\|_r^p\Bigg)^{1/p},$$
where, for all $\lambda \in \Lambda$, $\beta_{\lambda}$ stands for the column vector of $\mathbb{R}^r$ whose $l$-th line is $\beta_{\lambda}^{(l)}=\langle t^{(l)},\phi_{\lambda} \rangle_n$, for $l = 1,\ldots, r.$
\end{prop}

            \subsection{Proofs of Propositions~\ref{changenorm} and~\ref{maxf}}\label{sec:Proof23}

We present in a same section the proofs of Propositions~\ref{changenorm} and~\ref{maxf}, that both mainly call for the notion of decreasing rearrangement of a vector in $\mathbb{R}^n$.
\begin{defi}
Let $u\in\mathbb{R}^n$. The decreasing rearrangement of $u$ is the  $\mathbb{R}^n$- vector denoted by $u^\star$ satisfying
$$u^\star_1 \geq u^\star_2 \geq \ldots \geq u^\star_n \:\text{   and    }\: \{u^\star_i; 1\leq i \leq n\} = \{|u_i|; 1\leq i \leq n\}.$$
\end{defi}
\noindent
We will also make use of the Lorentz (quasi-)norms on $\mathbb{R}^n$ in the proof of Proposition~\ref{changenorm}, whose definition we recall here.
\begin{defi}
Let $0<p<+\infty$ and $0<q\leq +\infty$. We denote by $\|.\|_{\ell_{p,q}}$ the Lorentz (quasi-)norm defined on $\mathbb{R}^n$ by:
\begin{itemize}
\item [\textbullet] if $q$ is finite, $\|u\|_{\ell_{p,q}}=\Big(\sum_{i=1}^n i^{-1}(i^{1/p}u_i^\star)^q\Big)^{1/q}$;
\item [\textbullet] if $q=+\infty$, $\|u\|_{\ell_{p,\infty}}=\sup_{1\leq i \leq n}i^{1/p}u_i^\star$.
\end{itemize}
\end{defi}
\noindent
For all subset $I$ of $\{1,\ldots,n\}$, we denote by $\|.\|_{\ell_{p,q},I}$ the restriction of $\|.\|_{\ell_{p,q}}$ to $I$.
In particular, notice that, for all $u\in\mathbb{R}^n$, $0<p<+\infty$ and $0<q\leq +\infty$,
$$\|u\|_{\ell_{p,p}}=\|u\|_{\ell_p} \text{ and }
\|u^\star\|_{\ell_{p,q}}=\|u\|_{\ell_{p,q}}.$$
The reader may find in the appendix other useful properties relative to these notions.

            \subsubsection{Proof of Proposition~\ref{changenorm}}

The proof of Proposition~\ref{changenorm} mostly relies on a lemma that we demonstrate in this paragraph, after introducing a few notations.
Let $I$ be a dyadic interval of $\{1,\ldots,n\}$, $t\in\mathscr{M}(r,n)$, and $p>0$. By a compactness argument, there exists at least one vector in $\mathbb{R}^r$, denoted by $v_p(t,I)$, realizing the error $\mathcal{E}_p(t,I)$, i.e. satisfying
\begin{equation}
\mathcal{E}_p(t,I)=\bigg(\sum_{k\in I}\|t_k-v_p(t,I)\|_r^p\bigg)^{1/p}.
\notag
\end{equation}
We define the vectors $u_p(t,I)$ and $t^{\sharp,\alpha,p,I}$ in $\mathbb{R}^n$ whose coordinates are null outside of $I$ and given otherwise respectively by
\begin{equation}
\big(u_p(t,I)\big)_i=\|t_i-v_p(t,I)\|_r,\: \text{ for } i\in I,
\notag
\end{equation}
and
\begin{equation}
t^{\sharp,\alpha,p,I}_i = \sup_{I \supset J \ni i}|J|^{-(\alpha+1/p)}\mathcal{E}_p(t,J), \text{ for } i\in I,
\notag
\end{equation}
where the supremum is taken over all the dyadic intervals $J$ of $\{1,\ldots,n\}$ that are contained in $I$ and contain $i$.
\begin{lemm}\label{changenormlemm}
Let $\alpha > 0$, $p>0$ and $t\in \mathscr{M}(r,n)$.
Let $I$ be a dyadic interval of $\{1,\ldots,n\}$ containing at least two elements.
For all $j\in \{1,\ldots,|I|/2\}$,
\begin{equation}
\big(u_p(t,I)\big)^\star_j \leq C(\alpha,p) \Bigg(\sum_{k=j}^{|I|/2}k^{\alpha-1}\big(t^{\sharp,\alpha,p,I}\big)^\star_k +j^\alpha\big(t^{\sharp,\alpha,p,I}\big)^\star_j\Bigg).
\label{eq:ineqchangenormlemm}
\notag
\end{equation}
\end{lemm}
\begin{proof}
We fix $j\in \{1,\ldots,|I|/2\}$.
Let $E$ be the set composed of all the indices $i$ in $\{1,\ldots,n\}$ satisfying
$(t^{\sharp,\alpha,p,I})_i > (t^{\sharp,\alpha,p,I})_{j}^\star$. As $|E|\leq j-1$, we only have to prove that
\begin{equation}
\big(u_p(t,I)\big)_i \leq C(\alpha,p) \Bigg(\sum_{k=j}^{|I|/2}k^{\alpha-1}\big(t^{\sharp,\alpha,p,I}\big)^\star_k +j^\alpha\big(t^{\sharp,\alpha,p,I}\big)^\star_{j}\Bigg)
\label{eq:ineq2changenormlemm}
\end{equation}
for all the indices $i\in\{1,\ldots,n\}$, except maybe for those belonging to $E$.
Consider $i \in \{1,\ldots,n\}$ such that $i \notin E$. If $i\notin I$, then $\big(u_p(t,I)\big)_i=0$, so Inequality~\eqref{eq:ineq2changenormlemm} is trivial. Suppose now that $i\in I$ and $i\notin E$, and let $\{I_l\}_{1\leq l\leq m}$ be the sequence of dyadic intervals defined by
$$I_1=I,\: I_{l+1} \text{ is the child of } I_l \text{ containing } i, \text { and } I_m=\{i\},$$
where $m\geq 2$ because $|I|\geq 2$.
Notice that, for all $l\in\{0,\ldots,m-1\}$, $|I_{l+1}|=2^{-l}|I|$.
Let $q$ be the strictly positive integer such that
$$2^{-(q+1)}|I| < j\leq 2^{-q}|I|.$$
Such a definition implies, in particular, that $2^{-q}|I|\geq 1$, so that $q<m$.
From the triangular inequality,
\begin{equation}
\big(u_p(t,I)\big)_i
\leq
\sum_{l=2}^q \|v_p(t,I_{l-1})-v_p(t,I_l)\|_r
+ \sum_{l=q+1}^m \|v_p(t,I_{l-1})-v_p(t,I_l)\|_r,
\label{eq:triangineq2changenormlemm}
\end{equation}
with the convention that the first sum in Inequality~\eqref{eq:triangineq2changenormlemm} is null for $q=1$.
Let us fix $l\in\{2,\ldots,m\}$ and determine an upper-bound for the term $\|v_p(t,I_{l-1})-v_p(t,I_l)\|_r$. We recall  that $I_l\subset I_{l-1}$ and $|I_{l-1}|=2|I_l|$. Besides, for all $p>0$, the (quasi-)norm $\|.\|_{\ell_p}$ satisfies a triangular inequality within a multiplicative constant $C(p)$, where we can take $C(p)=1$ for $p\geq1$, and $C(p)=2^{1/p}$ for $0<p<1$. Therefore, we get
\begin{equation}
\|v_p(t,I_{l-1})-v_p(t,I_l)\|_r
\leq C(p)|I_l|^{-1/p}\Big(\mathcal{E}_p(t,I_{l-1})+\mathcal{E}_p(t,I_l)\Big),
\notag
\end{equation}
which leads to
\begin{equation}
\|v_p(t,I_{l-1})-v_p(t,I_l)\|_r \leq C(\alpha,p)|I_l|^{\alpha} \min_{k\in I_l} t^{\sharp,\alpha,p,I}_k.
\label{eq:ineq3changenormlemm}
\end{equation}
Let us bound the first sum appearing in~\eqref{eq:triangineq2changenormlemm}.
For all $l \in \{2,\ldots,m\}$, we have
\begin{equation}
\min_{k\in I_l} t^{\sharp,\alpha,p,I}_k
\leq \big(t^{\sharp,\alpha,p,I}\big)^\star_{|I_l|}
= \min_{1\leq k\leq |I_l|}\big(t^{\sharp,\alpha,p,I}\big)^\star_k,
\notag
\end{equation}
and, as $|I_{l+1}|=|I_l|/2$,
\begin{equation}
|I_l|^\alpha
= C(\alpha)\int_{|I_{l+1}|}^{|I_l|}x^{\alpha-1}\,\mathrm{d}x
\leq C(\alpha) \sum_{k=|I_{l+1}|}^{|I_l|}k^{\alpha-1}.
\notag
\end{equation}
Consequently, when $q\geq2$, Inequality~\eqref{eq:ineq3changenormlemm} yields
\begin{align}
\sum_{l=2}^q \|v_p(t,I_{l-1})-v_p(t,I_l)\|_r
&\leq C(\alpha,p) \sum_{l=2}^q \sum_{k=|I_{l+1}|}^{|I_l|}k^{\alpha-1}\big(t^{\sharp,\alpha,p,I}\big)^\star_k \notag\\
&\leq C(\alpha,p) \sum_{k=j}^{|I|/2}k^{\alpha-1}\big(t^{\sharp,\alpha,p,I}\big)^\star_k.
\notag
\end{align}
Regarding the second sum appearing in~\eqref{eq:triangineq2changenormlemm}, we now use Inequality~\eqref{eq:ineq3changenormlemm} combined with the following remarks. For all $l$ such that $q+1\leq l\leq m$, we have $\min_{k\in I_l} t^{\sharp,\alpha,p,I}_k \leq t^{\sharp,\alpha,p,I}_i$, since $I_l$ contains $i$, and we recall that $|I_l|=2^{-(l-1)}|I|$. Therefore,
\begin{equation}
\sum_{l=q+1}^m \|v_p(t,I_{l-1})-v_p(t,I_l)\|_r
\leq C(\alpha,p)|I|^\alpha \big(t^{\sharp,\alpha,p,I}\big)_i
\sum_{l=q+1}^m 2^{-(l-1)\alpha}.
\notag
\end{equation}
Furthermore, remember that $2^{-(q+1)}|I|<j$ and $i \notin E$, so we finally obtain
\begin{equation}
\sum_{l=q+1}^m \|v_p(t,I_{l-1})-v_p(t,I_l)\|_r
\leq C(\alpha,p) j^\alpha\big(t^{\sharp,\alpha,p,I}\big)_{j}^\star.
\notag
\end{equation}
We have thus proved inequality~\eqref{eq:ineq2changenormlemm} and Lemma~\ref{changenormlemm}.
\end{proof}

We are now able to prove Proposition~\ref{changenorm}.
Let $\alpha>0$, $p(\alpha)=(\alpha+1/2)^{-1}$ and $t\in\mathscr{M}(r,n)$. We fix $i\in\{1,\ldots,n\}$. From the definition of $\mathcal{E}_2(t,I)$ for any subset $I$ of $\{1,\ldots,n\}$, and due to the fact that $\mathcal{E}_2(t,\{i\})=0$, we have
$$t^{\sharp,\alpha,2}_i \leq \sup_{I\ni i}|I|^{-1/p(\alpha)}\|u_{p(\alpha)}(t,I)\|_{\ell_2},$$
where the supremum is taken over all the dyadic intervals $I$ of $\{1,\ldots,n\}$ that contain $i$, except for $\{i\}$.
We fix such an interval $I$.
The sequence $\big\{\big(u_{p(\alpha)}(t,I)\big)^\star_j\big\}_{1\leq j \leq n}$ decreases and is null for $j\geq |I|+1$, hence
\begin{equation}
\big\|u_{p(\alpha)}(t,I)\big\|_{\ell_2}^2
\leq 2\sum_{j=1}^{|I|/2} \Big(\big(u_{p(\alpha)}(t,I)\big)^\star_j\Big)^2.
\notag
\end{equation}
From Lemma~\ref{changenormlemm} and the definition of $p(\alpha)$, we get
\begin{equation}
\big\|u_{p(\alpha)}(t,I)\big\|_{\ell_2}^2
\leq C(\alpha)
\Bigg(\sum_{j=1}^{|I|/2}j^{-1}\Bigg(j^{1/2}\sum_{k=j}^{|I|/2}k^{\alpha-1}\big(t^{\sharp,\alpha,p(\alpha),I}\big)^\star_k\Bigg)^2 +
\big\|\big(t^{\sharp,\alpha,p(\alpha),I}\big)^\star\big\|_{\ell_{p(\alpha),2}}^2\Bigg).
\notag
\end{equation}
Using one of Hardy's inequalities (cf. Proposition~\ref{HardyProp} in the Appendix) and noticing that $t^{\sharp,\alpha,p(\alpha),I}\leq t^{\sharp,\alpha,p(\alpha)}$, we are led to
\begin{equation}
\big\|u_{p(\alpha)}(t,I)\big\|_{\ell_2}
\leq C(\alpha) \|t^{\sharp,\alpha,p(\alpha)}\|_{\ell_{p(\alpha),2},I}.
\notag
\end{equation}
Last, since $p(\alpha)<2$, we deduce from classical inequalities between Lorentz (quasi-)norms (cf. Proposition~\ref{Lorentzineq} in the Appendix)
$$t^{\sharp,\alpha,2}_i \leq C(\alpha)\sup_{I\ni i}|I|^{-1/p(\alpha)}\|t^{\sharp,\alpha,p(\alpha)}\|_{\ell_{p(\alpha)},I}$$
where the supremum is taken over all the dyadic intervals $I$ of $\{1,\ldots,n\}$ that contain $i$, which completes the proof of Proposition~\ref{changenorm}.

                        \subsubsection{Proof of Proposition~\ref{maxf}}

Let $q>1$ and $u\in \mathbb{R}^n$. As $M_1(u)=M_1(|u|)$, we can suppose that $u$ has positive or null coordinates.
Let us first demonstrate that, for all $i \in \{1,\ldots,n\}$,
\begin{equation}
(M_1(u))^{\star}_i \leq C \Bigg(i^{-1}\sum_{k=1}^i u^{\star}_k\Bigg).
\label{eq:decrearrmaxfunc}
\end{equation}
If $i=1$, then this inequality easily follows from the definitions of $(M_1(u))^{\star}_1$ and $u^{\star}_1.$ Let us now fix $i \in \{2,\ldots,n\}$.
We can write $u$ as $u=v+w$, where $v$ and $w$ are the $\mathbb{R}^n$-vectors whose respective coordinates are
$$v_k=\max\{u_k-u_{i}^\star,0\} \text{   and   } w_k=\min\{u_k,u_{i}^\star\}, \text{ for } k = 1,\ldots,n.$$
From the triangular inequality, we deduce that $M_1(u)\leq M_1(v)+M_1(w)$. Proposition~\ref{decrearrprop} (cf. Appendix) then leads to
\begin{equation}
(M_1(u))_{i}^\star \leq (M_1(v))^\star_{\lceil i/2\rceil}+(M_1(w))^\star_{\lfloor i/2\rfloor}.
\label{eq:maxflemm1int3}
\notag
\end{equation}
Moreover,
$$(M_1(w))^\star_{\lfloor i/2\rfloor} \leq \|M_1(w)\|_{\ell_{\infty}} \leq \|w\|_{\ell_{\infty}},$$
and, from Proposition~\ref{decrearrprop} again,
$$(M_1(v))^\star_{\lceil i/2 \rceil}\leq 2{i}^{-1}\|v\|_{\ell^1}.$$
Consequently,
\begin{equation}
(M_1(u))^{\star}_{i} \leq C\big({i}^{-1}\|v\|_{\ell^1}+\|w\|_{\ell^{\infty}}\big).
\label{eq:maxflemm1}
\end{equation}
Let $I$ be the set of all the indices $l$, $1\leq l\leq n$, such that $u_l > u_{i}^\star$. From the definitions of $v$ and $w$, we get
\begin{align}
\|v\|_{\ell^1}+i\|w\|_{\ell^{\infty}}
&\leq \sum_{k=1}^{|I|}u_k^\star+(i-|I|)u_{i}^\star=\sum_{k=1}^{i}u_k^\star,
\notag
\end{align}
which, given Inequality~\eqref{eq:maxflemm1}, completes the proof of~\eqref{eq:decrearrmaxfunc}.
We now have
\begin{equation}
\|(M_1(u))^{\star}\|_{\ell_q}^q
\leq C(q)\sum_{i=1}^n\Bigg(i^{-1}\sum_{k=1}^iu_k^\star\Bigg)^q.
\label{eq:maxfintineq}
\end{equation}
Let us denote by $q'$ the conjugate exponent of $q$, and write, for all $k$ in $\{1,\ldots,n\}$, $u_k^\star=k^{-1/qq'}k^{1/qq'}u_k^\star$. We deduce from Hölder's inequality
\begin{equation}
\sum_{i=1}^n\Bigg(i^{-1}\sum_{k=1}^iu_k^\star\Bigg)^q
\leq\sum_{i=1}^n\Bigg(q'i^{-1/q}\Bigg)^{q/q'}
\Bigg(i^{-1}\sum_{k=1}^ik^{1/q'}{(u_k^\star)}^q\Bigg).
\notag
\end{equation}
Interchanging the order of the summations, we obtain
\begin{equation}
\sum_{i=1}^n\Bigg(i^{-1}\sum_{k=1}^iu_k^\star\Bigg)^q
\leq C(q)\sum_{k=1}^n {(u_k^\star)}^q.
\notag
\end{equation}
Consequently,
$$\|(M_1(u))^{\star}\|_{\ell_q} \leq C(q)\|u^\star\|_{\ell_q},$$
hence Proposition~\ref{maxf}.

                \subsection{Proof of Proposition~\ref{compSharpBesov}}\label{sec:Proof4}

Let $p\in (0,2]$, $\alpha > 1/p-1/2$ and $t \in \mathscr{M}(r,n)$.
For all $i \in \{1,\ldots,n\}$ and all $0\leq J\leq N$, we denote by $I(J,i)$ the only dyadic interval of length $n2^{-J}$ that is contained in $\{1,\ldots,n\}$ and contains $i$. From the definition of $t^{\sharp,\alpha,p}$, we deduce
\begin{equation}
\|t^{\sharp,\alpha,p}\|_{\ell_{p}}^p
\leq \sum_{J = 0}^{N-1} (n^{-1}2^J)^{\alpha p + 1}\sum_{i=1}^n\Big(\mathcal{E}_p\big(t,I(J,i)\big)\Big)^p.
\label{eq:tsharplpnorm}
\end{equation}

Let us first suppose that $0<p\leq1$. From the definition of $\mathcal{E}_p(t,I(J,i))$, we have
\begin{equation}
\Big(\mathcal{E}_p\big(t,I(J,i)\big)\Big)^p
\leq \sum_{k \in I(J,i)}\|t_k-t_i\|_r^p.
\notag
\end{equation}
For all $-1 \leq j \leq N-1$, the functions $\{\phi_\lambda\}_{\lambda \in \Lambda(j)}$ are constant over any dyadic interval of length $n2^{-(j+1)}$. Therefore, if $k$ belongs to $I(J,i)$, then
\begin{equation}
t_k-t_i = \sum_{j=J}^{N-1}\sum_{\lambda \in \Lambda(j)}\beta_\lambda(\phi_{\lambda \: k} - \phi_{\lambda \: {i}}).
\notag
\end{equation}
As $0<p\leq 1$, we deduce from the classical inequality between $\ell_p$-quasi-norm and $\ell_1$-norm
\begin{align}
\sum_{i=1}^n\Big(\mathcal{E}_p\big(t,I(J,i)\big)\Big)^p
&\leq 2n^{2-p/2}2^{-J}\sum_{j=J}^{N-1}2^{jp(1/2-1/p)}\sum_{\lambda \in \Lambda(j)}\|\beta_\lambda\|_r^p.
\notag
\end{align}
Interchanging the order of the summations, we get
\begin{align}
\|t^{\sharp,\alpha,p}\|_{\ell_{p}}^p
&\leq C(\alpha,p)n^{1-p(\alpha +1/2)}\sum_{j = 0}^{N-1}2^{jp(\alpha + 1/2 - 1/p)}\sum_{\lambda \in \Lambda(j)}\|\beta_\lambda\|_r^p.
\notag
\end{align}

Let us now consider the case $1<p\leq 2$. We fix $0\leq J\leq N-1$ and define
$$T(J)=\sum_{j=J}^{N-1}\sum_{\lambda \in \Lambda(j)}\beta_\lambda\phi_{\lambda}.$$
As $t-T(J)$ is constant over any dyadic interval of length $n2^{-J}$,
\begin{equation}
\mathcal{E}_p\big(t,I(J,i)\big)
=\mathcal{E}_p\big(T(J),I(J,i)\big).
\notag
\end{equation}
This equality and the definition of $\mathcal{E}_p\big(T(J),I(J,i)\big)$ lead to
\begin{align}
\sum_{i=1}^n\Big(\mathcal{E}_p\big(t,I(J,i)\big)\Big)^p
&\leq \sum_{i=1}^n\sum_{k\in I(J,i)}\big\|\big(T(J)\big)_k\big\|_r^p\notag\\
&\leq n2^{-J}\sum_{k=1}^n\bigg(\sum_{j=J}^{N-1}\sum_{\lambda \in \Lambda(j)}\|\beta_\lambda\|_r|\phi_{\lambda\:k}|\bigg)^p.
\notag
\end{align}
From~\eqref{eq:tsharplpnorm} and this last inequality, we get
\begin{equation}
\|t^{\sharp,\alpha,p}\|_{\ell_{p}}^p
\leq n^{-\alpha p}\sum_{k=1}^n\sum_{J = 0}^{N-1} \bigg(2^{J\alpha }\sum_{j=J}^{N-1}\sum_{\lambda \in \Lambda(j)}\|\beta_\lambda\|_r|\phi_{\lambda\:k}|\bigg)^p.
\notag
\end{equation}
Then, using one of Hardy's inequalities (cf. Proposition~\ref{HardyProp} in the Appendix) and remembering that, for all $j\in\{-1,\ldots,N-1\}$, the functions $\{\phi_\lambda\}_{\lambda\in\Lambda(j)}$ have disjoint supports, we conclude that
\begin{align}
\|t^{\sharp,\alpha,p}\|_{\ell_{p}}^p
&\leq C(\alpha,p)n^{-\alpha p}\sum_{j = 0}^{N-1}2^{j\alpha p}\sum_{\lambda \in \Lambda(j)}\|\beta_\lambda\|_r^p\sum_{k=1}^n|\phi_{\lambda\:k}|^p,
\notag
\end{align}
hence Proposition~\ref{compSharpBesov}.
\clearpage

\appendix

\section{Some useful inequalities}

We state here, for vectors in $\mathbb{R}^n$, a few inequalities that are similar to classical inequalities for functions of a continuous parameter. The proofs of the latter, which may be found in~\cite{BennShar}, for instance, are easy to transpose to the finite-dimensional case.

\begin{prop}[Some properties of decreasing rearrangements]\label{decrearrprop}
Let $u$ and $v$ be two vectors in $\mathbb{R}^n$.
For all $\lambda \geq 0$, let $I_u(\lambda)$ be the set of the indices $k$ in $\{1,\ldots,n\}$ such that $|u_k| \geq \lambda$.
\begin{enumerate}
\item For all $i \in \{1,\ldots,n\}$,
$u^\star_i=\sup\{\lambda\geq 0 \:s.t.\: |I_u(\lambda)| \geq i\}.$
\item If, for all $i\in \{1,\ldots,n\}$, $u_i\leq v_i$, then, for all $i\in \{1,\ldots,n\}$, $u_i^\star\leq v_i^\star.$
\item For all $i,j \in \{1,\ldots,n\}$ such that $1\leq i+j\leq n$, $(u+v)^\star_{i+j}\leq u_i^\star+v_j^\star.$
\item For all $i\in \{1,\ldots,n\}$, $(M_1(u))^\star_i \leq i^{-1}\|u\|_{\ell_1}.$
\end{enumerate}
\end{prop}
\begin{proof}
See, for instance,~\cite{BennShar}, Proposition 1.7. and Theorem 3.3.
\end{proof}

\begin{prop}[Inequalities between Lorentz (quasi-)norms]\label{Lorentzineq}
Let $p, q$ and $q'$ be positive reals and let $u$ be a vector in $\mathbb{R}^n$.
\begin{enumerate}
\item If $p\leq q$, then $\|u\|_{\ell_{p,\infty}} \leq C(p,q) \|u\|_{\ell_{p,q}}$.
\item If $q'\leq q$, then $\|u\|_{\ell_{p,q}} \leq C(p,q,q') \|u\|_{\ell_{p,q'}}.$
\end{enumerate}
\end{prop}
\begin{proof}
See, for instance,~\cite{BennShar}, Proposition 4.2.
\end{proof}

\begin{prop}[Hardy's inequalities]\label{HardyProp}
Let $q >1$ and let $\psi$ be a vector in $\mathbb{R}^n$ whose coordinates are non-negative.
\begin{enumerate}
\item For all $\lambda < 1$,
\begin{equation}
\sum_{i=1}^n i^{-1}\Bigg(i^{1-\lambda}\sum_{k=i}^n k^{-1}\psi_k\Bigg)^q
\leq C(\lambda,q) \sum_{i=1}^n i^{-1}(i^{1-\lambda}\psi_i)^q.
\label{eq:Hardy1}
\notag
\end{equation}
\item For all $\alpha>0$,
\begin{equation}
\sum_{i=1}^n \Bigg(2^{i\alpha}\sum_{k=i}^n \psi_k\Bigg)^q
\leq C(\alpha,q) \sum_{i=1}^n (2^{i\alpha}\psi_i)^q.
\label{eq:Hardy2}
\notag
\end{equation}
\end{enumerate}
\end{prop}
\begin{proof}
See, for instance,~\cite{BennShar}, Lemma 3.9.
\end{proof}

\clearpage


\end{document}